\documentclass[11pt,a4paper]{article}

\usepackage{fullpage} 
\usepackage{parskip} 
\usepackage{tikz} 
\usepackage{amssymb,amsmath}
\usepackage{graphicx,color}
\usepackage{float,subfigure}
\usepackage[colorlinks=true]{hyperref}
\usepackage{mathptmx} 

\newtheorem{thm}{Theorem}[section]
\newtheorem{prop}[thm]{Proposition}
\newtheorem{lem}[thm]{Lemma}

\newtheorem{rem}[thm]{Remark}

\numberwithin{equation}{section}

\newcommand{\qqed}{\hfill $\square$} 
\newcommand{\rom}[1]{\uppercase\expandafter{\romannumeral #1\relax}} 

\newcommand{\sgn}{\operatorname{sgn}}
\title{Isoperimetric deformations of curves on the Minkowski plane}

\makeatletter

\renewcommand\@date{{%
  \vspace{-\baselineskip}%
  \large\centering

  \begin{tabular}{@{}c@{}}
     Hyeongki Park \\[0.3pc]
    \normalsize Graduate school of Mathematics, Kyushu University \\
    \normalsize 744 Motooka, Fukuoka 819-0395, Japan \\
    \normalsize {\it h-paku@math.kyushu-u.ac.jp }
  \end{tabular}
  \\[1pc]

  \begin{tabular}{@{}c@{}}
     Jun-ichi Inoguchi \\[0.3pc]
    \normalsize Institute of Mathematics, University of Tsukuba \\
    \normalsize Tsukuba 305-8571, Japan \\
    \normalsize {\it inoguchi@math.tsukuba.ac.jp }
  \end{tabular}
  \\[1pc]

  \begin{tabular}{@{}c@{}}
     Kenji Kajiwara \\[0.3pc]
    \normalsize Institute of Mathematics for Industry, Kyushu University \\
    \normalsize 744 Motooka, Fukuoka 819-0395, Japan \\
    \normalsize {\it kaji@imi.kyushu-u.ac.jp }
  \end{tabular}
  \\[1pc]

  \begin{tabular}{@{}c@{}}
     Ken-ichi Maruno \\[0.3pc]
    \normalsize Department of Applied Mathematics School of Fundamental Science and Engineering,\\
    \normalsize Faculty of Science and Engineering Waseda University \\
    \normalsize 3-4-1 Okubo, Tokyo 169-8555, Japan \\
    \normalsize {\it kmaruno@waseda.jp }
  \end{tabular}
  \\[1pc]

  \begin{tabular}{@{}c@{}}
     Nozomu Matsuura \\[0.3pc]
    \normalsize Department of Education and Creation Engineering, Kurume Institute of Technology \\
    \normalsize 2228-66 Kamitsu-machi, Fukuoka 830-0052, Japan \\
    \normalsize {\it nozomu@kurume-it.ac.jp }
  \end{tabular}
  \\[1pc]

  \begin{tabular}{@{}c@{}}
     Yasuhiro Ohta \\[0.3pc]
    \normalsize Department of Mathematics, Kobe University \\
    \normalsize Rokko, Kobe 657-8501, Japan \\
    \normalsize {\it ohta@math.sci.kobe-u.ac.jp }
  \end{tabular}
  \\[1pc]

}}

\makeatother

\begin{document}

\maketitle

\begin{abstract}
  We formulate an isoperimetric deformation of curves on the Minkowski plane,
  which is governed by the defocusing mKdV equation.
  Two classes of exact solutions
  to the defocusing mKdV equation
  are also presented in terms of the $\tau$ functions.
  By using one of these classes,
  we construct an explicit formula for the corresponding motion of curves on the Minkowski plane
  even though those solutions have singular points.
  Another class give regular solutions to the defocusing mKdV equation.
  Some pictures illustrating typical dynamics of the curves are presented.
\end{abstract}

\section{Introduction}

{It is well-known that a certain class of integrable systems describes motions
of plane and space curves in various settings.}  For {instance}, the nonlinear
Schr\"{o}dinger equation describes {a motion of space curves which is a physical
model of vortex filaments \cite{hasimoto:1972}}, and the {modified Korteweg-de Vries
(mKdV) equation describes motion of space and plane curves preserving the arc length}
\cite{goldstein:1991,lamb:1976}.  {Curve motions have been studied not only in the Euclidean geometry but also in various Klein
geometries}
\cite{chou:2002,chou:2003,chou:2004,fujioka:2010,kajiwara:2016}.
For example, the KdV equation describes a motion of plane curves that
preserves the areal velocity \cite{pinkall:1995}.  Moreover, in recent years, explicit formulas for
these curve motions have been established in \cite{hirose:,inoguchi:2012,inoguchi:2012:2,park:2018} by using the
theory of $\tau$ functions.

In this paper, we consider the isoperimetric motions of curves on the Minkowski plane,
namely, motions preserving the arc length,
and show that the simplest nontrivial motion is described by the
{\em defocusing} mKdV equation,
\begin{equation}\label{eqn:defocusing_mKdV:0}
 \frac{\partial \kappa}{\partial t} = \frac{\partial^3 \kappa}{\partial x^3} - \frac{3}{2}\kappa^2\frac{\partial\kappa}{\partial x}.
\end{equation}
The regular soliton type solutions of the defocusing mKdV equation show rather different behavior
from other equations; it admits the solutions where the solitons run on a shock wave
\cite{perelman:1974,perelman:1974:2}. On the other hand, the {\em focusing} mKdV equation
\begin{equation}\label{eqn:focusing_mKdV:0}
 \frac{\partial \hat\kappa}{\partial t} = \frac{\partial^3 \hat\kappa}{\partial x^3}
+ \frac{3}{2}\hat\kappa^2\frac{\partial\hat\kappa}{\partial x},
\end{equation}
which governs the curve motions preserving the arc length in the Euclidean plane describes the
ordinary dynamics of solitons.  Those two mKdV equations cannot be transformed to each other by
scale change. It should be remarked that the defocusing mKdV equation has not been studied well
because of less physical relevance compared to the focusing mKdV equation.

This paper is organized as follows.
In Section \ref{section:curve},
we review some basic notions of curves on the Minkowski plane.
In Section \ref{section:deformation},
we formulate an integrable isoperimetric deformation of curves on the Minkowski plane.
In Section \ref{section:solution1},
we first construct a class of exact solutions to the defocusing mKdV equation
in terms of the $\tau$ functions.
Those solutions, however, give singular solutions.
We also construct regular solutions to the defocusing mKdV equation
in Section \ref{section:solution2}
by suitable choice of parameters and reductions on the $\tau$ functions, {following the
idea given in \cite{ohta:1989}.}
We finally present the formulas for corresponding curve motions on the Minkowski plane.

\section{Curves on the Minkowski plane}\label{section:curve}

In this section,
we state a curve theory on the Minkowski plane.
We provide the vector space $\mathbb{R}^2$ with
the Lorentzian inner product $\left\langle\,\cdot\,,\,\cdot\,\right\rangle$
\begin{equation*}
\left\langle v, w\right\rangle
= v_1 w_1 - v_2 w_2
\end{equation*}
for arbitrary vectors $v = {}^\mathrm{t} \left[v_1, v_2\right]$
and $w = {}^\mathrm{t} \left[w_1, w_2\right]$.
We write $\mathbb{R}^{1,1}$ for
$\left(\mathbb{R}^2, \left\langle\,\cdot\,,\,\cdot\,\right\rangle\right)$,
and call it the \textit{Minkowski plane}.
We say that a vector $v \in \mathbb{R}^{1,1}$ is
\textit{spacelike} if $\langle v, v\rangle > 0$ or $v=0$,
\textit{timelike} if $\langle v, v\rangle < 0$,
and \textit{lightlike} if $\langle v, v\rangle =0$ and $v\neq0$.
For a vector $v \in \mathbb{R}^{1,1}$,
its norm $\left|v\right|$ is defined as
\begin{equation*}
\left|v\right| =
\begin{cases}
\left\langle v, v\right\rangle^{1/2} & \text{if $v$ is spacelike},\\
\left(- \left\langle v, v\right\rangle\right)^{1/2} &
\text{if $v$ is timelike},\\
0 & \text{if $v$ is lightlike}.
\end{cases}
\end{equation*}
Therefore
a spacelike (resp.\ timelike) vector
$v \in \mathbb{R}^{1,1}$ is of unit length
if and only if
$v \in H^1_+$
(resp.\ $H^1_-$),
where $H^1_\pm$ are the hyperbolas of two sheets
\begin{align*}
H^1_+ =
\left\{\left.
v \in \mathbb{R}^{1,1}
\;\right|\;
\left\langle v,v\right\rangle = 1\right\},\quad
H^1_- =
\left\{\left.
v \in \mathbb{R}^{1,1}
\;\right|\;
\left\langle v,v\right\rangle = -1\right\}.
\end{align*}
The \textit{Lorentz group}
\begin{equation*}
\mathrm{O} \left(1,1\right)=\left\{
A\in \mathrm{GL}\left(2,\mathbb{R}\right)
\;\left|\; {}^\mathrm{t}\!A E' A = E'\right.\right\},\quad
E'=
\begin{bmatrix}
1 & 0\\
0 & -1
\end{bmatrix},
\end{equation*}
preserves the inner product.
Indeed the equality $\left\langle Av, Aw\right\rangle =
\left\langle v, w\right\rangle$ holds for all $v, w \in \mathbb{R}^{1,1}$
if and only if $A \in \mathrm{O}\left(1,1\right)$.
We also consider the subgroup $\mathrm{SO}\left(1,1\right)
= \mathrm{O}\left(1,1\right) \cap \mathrm{SL}\left(2,\mathbb{R}\right)$,
which has two connected components
\begin{align*}
\mathrm{SO}^+\left(1,1\right)
= \left\{\left.%
\begin{bmatrix}
a & b\\
c & d
\end{bmatrix}
\in \mathrm{SO}\left(1,1\right) \;\right|\; d > 0\right\}
&= \left\{\left.%
\begin{bmatrix}
\cosh \phi & \sinh \phi\\
\sinh \phi & \cosh \phi
\end{bmatrix}
\;\right|\; \phi \in \mathbb{R}\right\},\\
\mathrm{SO}^-\left(1,1\right)
= \left\{\left.%
\begin{bmatrix}
a & b\\
c & d
\end{bmatrix}
\in \mathrm{SO}\left(1,1\right) \;\right|\; d < 0\right\}
&= \left\{\left.%
\begin{bmatrix}
- \cosh \phi & \sinh \phi\\
\sinh \phi & - \cosh \phi
\end{bmatrix}
\;\right|\; \phi \in \mathbb{R}\right\}.
\end{align*}

Let $\gamma\colon I\subset\mathbb{R}\rightarrow \mathbb{R}^{1,1},\,
\xi \mapsto \gamma \left(\xi\right)$ be a regular curve
on the Minkowski plane,
parametrized by an arbitrary parameter $\xi$.
We say that $\gamma$ is \textit{spacelike} if
its velocity $\left(d/d\xi\right)\gamma$ is spacelike everywhere.
A spacelike curve $\gamma$ is said to be \textit{unit-speed}
if its velocity is of unit length everywhere.
Therefore,
the velocity $T$ of a unit-speed spacelike curve $\gamma$
moves along the hyperbola of two sheets $H^1_+$.
Since $T$ is a continuous vector field along $\gamma$,
$T$ moves along one of the sheets of $H^1_+$.
Similarly the notion of unit-speed timelike curve is defined.
Hereafter,
we consider a unit-speed spacelike curve
\begin{equation*}
\gamma \left(x\right) =
\begin{bmatrix}
\gamma_1 \left(x\right)\\
\gamma_2 \left(x\right)
\end{bmatrix},
\end{equation*}
and assume that ${\gamma_1}' \left(x\right) > 0$.
We say that such a tangent vector field $T = \gamma'$ is
\textit{positive pointing}.
The positive pointing tangent vector field $T$ is obviously expressed as
\begin{equation}\label{positivepointing}
T=
\begin{bmatrix}
\cosh \theta\\
\sinh \theta
\end{bmatrix}
\end{equation}
with some function $\theta$.
We define the \textit{normal vector field} by
\begin{equation*}
N = J' T =
\begin{bmatrix}
\sinh \theta\\
\cosh \theta
\end{bmatrix},\quad
J' =
\begin{bmatrix}
0 & 1\\
1 & 0
\end{bmatrix}.
\end{equation*}
They satisfy
$\left|T\right| = \left|N\right| = 1$ and $\langle T,N\rangle=0$.
Introducing a frame
$\Phi\colon I\rightarrow\mathrm{SO}^+\left(1,1\right)$
by $\Phi=\left[T,N\right]$, we have the Frenet formula
\begin{equation}\label{frenet}
\Phi'=\Phi L,\quad
L=
\begin{bmatrix}
0 & \theta'\\
\theta' & 0
\end{bmatrix}.
\end{equation}
On the other hand,
there exists a function $\kappa$ called the \textit{curvature} of $\gamma$,
such that $T' = \kappa N$
because $\left\langle T, T'\right\rangle=0$.
Thus we have $\kappa = \theta'$,
from which $\theta$ is sometimes referred to
as the \textit{potential} function.
The discussion is summarized as
the fundamental theorem of plane curves as follows.
\begin{prop}\label{prop:representation}
For a given function $\kappa = \kappa \left(x\right)$,
there exists a unit-speed spacelike curve $\gamma = \gamma \left(x\right)$
on $\mathbb{R}^{1,1}$ with a positive pointing tangent vector field,
such that $\kappa$ is the curvature of $\gamma$.
In fact $\gamma$ is given by the integral
\begin{equation}\label{representation}
\gamma \left(x\right) = \int_{x_0}^x
\begin{bmatrix}
\cosh \theta\\
\sinh \theta
\end{bmatrix}
dx,\quad
\theta \left(x\right) = \int_{x_0}^x \kappa\, dx.
\end{equation}
Moreover,
if two unit-speed spacelike curves $\gamma$, $\overline{\gamma}$
have the same curvature,
then they differ only by a Lorentzian motion,
namely there exists a matrix $A \in \mathrm{SO}(1,1)$
and a vector $v \in \mathbb{R}^{1,1}$
such that $\overline{\gamma} \left(x\right) =
A \gamma \left(x\right) + v$.
\end{prop}
\begin{rem}
The arclength function of a curve $\gamma(\xi)$ on the Minkowski plane is defined by
\begin{equation}\label{def:arclength}
  x \left(\xi, t\right) =
  \int_{\xi_0}^\xi \left|\gamma_\xi\right| d\xi.
\end{equation}
so that the arclength parametrized curve $\gamma(x)$ is unit-speed.
\end{rem}
\begin{rem}
In terms of an arbitrary parameter $\xi$,
we can rewrite Proposition \ref{prop:representation} as follows.
Let $\tilde{\kappa} = \tilde{\kappa}\left(\xi\right)$ be a function
and $x = x\left(\xi\right)$ be a monotonously increasing function.
Then,
up to Lorentzian motions,
there uniquely exists a spacelike curve
$\tilde{\gamma} = \tilde{\gamma} \left(\xi\right)$ on $\mathbb{R}^{1,1}$
with a positive pointing tangent vector field,
such that $\tilde{\kappa}$ is the curvature of $\tilde{\gamma}$
and $x$ is the arclength of $\tilde{\gamma}$.
In fact $\tilde{\gamma}$ is given by the integral
\begin{equation}\label{representation-gen}
\tilde{\gamma} \left(\xi\right) = \int_{\xi_0}^\xi
\begin{bmatrix}
\cosh \tilde{\theta}\,\\
\sinh \tilde{\theta}
\end{bmatrix}
x' d\xi,\quad
\tilde{\theta} \left(\xi\right)
= \int_{\xi_0}^\xi \tilde{\kappa}x' d\xi.
\end{equation}
\end{rem}

\section{Deformation of curves on the Minkowski plane}\label{section:deformation}

We formulate an arclength preserving deformation of
a unit speed spacelike curve $\gamma = \gamma \left(x\right)$
on $\mathbb{R}^{1,1}$
with positive pointing tangent vector field,
and show that it can be governed by the defocusing mKdV equation.
Introducing a deformation parameter $t$,
we denote again by $\gamma = \gamma \left(\xi, t\right)$ the deformation,
where $\gamma \left(\xi, 0\right)$ is the initial
unit-speed curve $\gamma \left(x\right)$.
We decompose $\left(\partial/\partial t\right) \gamma$ in the form
\begin{equation}\label{expression-of-curvemotion}
\gamma_t = f T + g N,
\end{equation}
where $T= \left|\gamma_\xi\right|^{-1} \gamma_\xi$ and $N = J' T$.
Here the subscripts mean differentiation with respect to
the indicated variables.
\begin{prop}
Let $\gamma = \gamma \left(\xi, t\right)$ be
a family of spacelike curves on $\mathbb{R}^{1,1}$
with positive pointing tangent vector fields,
such that the initial curve $\gamma \left(\xi, 0\right)$ is unit speed.
Then the arclength is independent of $t$ if and only if
$\left\langle\gamma_{t \xi}, \gamma_\xi\right\rangle = 0$.
\end{prop}
\noindent{\it Proof.}
Differentiating \eqref{def:arclength} by $t$,
we have
\begin{equation*}
x_t
= \int_{\xi_0}^\xi
\frac{\left\langle\gamma_{t \xi}, \gamma_\xi\right\rangle}%
{\left|\gamma_\xi\right|}\, d\xi.
\end{equation*}
Therefore $x_t = 0$ for all $\xi$ if and only if
$\left\langle\gamma_{t \xi}, \gamma_\xi\right\rangle = 0$.
\qqed\\
Since we have from \eqref{representation-gen} that
$T_\xi = x_\xi \kappa N$,
where $\kappa$ is the curvature of $\gamma$ at each $t$,
it follows from the expression \eqref{expression-of-curvemotion} that
\begin{equation*}
\gamma_{t \xi} =
\left(f_\xi + x_\xi \kappa g\right) T
+ \left(g_\xi + x_\xi \kappa f\right) N.
\end{equation*}
Therefore the isoperimetric condition
$\left\langle\gamma_{t \xi}, \gamma_\xi\right\rangle = 0$
is equivalent to the equality
$f_\xi + x_\xi \kappa g=0$.
In the followings,
we consider an isoperimetric deformation of
a unit-speed spacelike curve,
and hence we can assume that $\xi$ itself is the arclength parameter.
Thus,
the isoperimetric condition becomes
\begin{equation}\label{condition:isoperimetric}
f_x + \kappa g = 0,
\end{equation}
and the frame $\Phi=\left[T,N\right]$ is deformed as
\begin{equation}\label{deformation-of-frame-0}
\Phi_t = \Phi M,\quad
M =
\begin{bmatrix}
0 & g_x+\kappa f\\
g_x+\kappa f & 0
\end{bmatrix}.
\end{equation}
Under the isoperimetric condition \eqref{condition:isoperimetric},
the compatibility condition between \eqref{frenet}
and \eqref{deformation-of-frame-0},
$L_t-M_x-LM+ML=0$,
is
\begin{equation*}
\kappa_t
= \left(g_x + \kappa f\right)_x
= \Omega g,
\end{equation*}
where $\Omega = \partial_x^2
- \kappa^2 - \kappa_x \partial_x^{-1} \left(\kappa\,\cdot\,\right)$
is the recursion operator of the defocusing mKdV hierarchy.
In view of this,
it is reasonable to choose $g$ as $g=\kappa_x$
and hence $f=-\kappa^2/2$.
Thus we have:
\begin{thm}[defocusing mKdV flow]
Let $\gamma = \gamma \left(x,t\right)$ be a family
of unit speed spacelike curves on $\mathbb{R}^{1,1}$
with positive pointing vector field,
and $\kappa$ the curvature of $\gamma$ at each $t$.
Then $\gamma$ is an arclength preserving deformation,
and it varies according to the formula
\begin{equation}\label{mkdvflow}
\gamma_t = - \frac{\kappa^2}{2} T + \kappa_x N
\end{equation}
if and only if $\kappa$ satisfies the defocusing mKdV equation
\begin{equation}\label{mkdv}
\kappa_t = \kappa_{xxx} - \frac{3}{2}\kappa^2\kappa_x.
\end{equation}
\end{thm}
\noindent{\it Proof.}
The frame is deformed as
\begin{equation}\label{deformation-of-frame}
\Phi_t=\Phi M,\quad
M=
\begin{bmatrix} 0 &
\kappa_{xx}-\frac{1}{2} \kappa^3\\
\kappa_{xx}-\frac{1}{2} \kappa^3 & 0
\end{bmatrix}.
\end{equation}
The compatibility condition between
\eqref{frenet} and \eqref{deformation-of-frame} is \eqref{mkdv}.
\qqed\\
We call the deformation \eqref{mkdvflow}
the \textit{defocusing mKdV flow}.
If $\gamma$ is a defocusing mKdV flow,
then the potential function $\theta$ satisfies
the \textit{potential defocusing mKdV equation}
\begin{equation}\label{pot_mKdV}
\theta_t = \theta_{xxx} - \frac{1}{2} {\theta_x}^3.
\end{equation}

\section{Solutions}

We construct solutions to the defocusing mKdV equation
in terms of $\tau$ functions,
and derive an explicit formula
for the defocusing mKdV flow \eqref{mkdvflow}.

\subsection{Explicit formula}\label{section:solution1}

Let $\tau=\tau\left(x,t;y\right)$ and
$\overline{\tau}=\overline{\tau}\left(x,t;y\right)$
be real-valued functions,
where $y$ is an auxiliary variable.
For a real constant $c$,
we consider the system of bilinear equations
\begin{align}
&D_{x}D_{y}\,\tau\cdot\tau=-2\,{\overline{\tau}}^2,\label{c_bi_1}\\
&D_{x}D_{y}\,\overline{\tau}\cdot\overline{\tau}=-2\,\tau^2,\label{c_bi_1-2}\\
&\left( D_x^2-c \right)\,\tau\cdot\overline{\tau}=0,\label{c_bi_2}\\
&\left( D_x^3-D_t-3cD_x \right) \tau\cdot\overline{\tau}=0.\label{c_bi_3}
\end{align}
Here $D_x$, $D_y$ and $D_t$ are
Hirota's bilinear differential operators
\cite{hirota:2004},
defined as
\begin{equation*}
D_x^{i}D_y^{j}\, f\cdot g
=\left.\left(\frac{\partial}{\partial x}
-\frac{\partial}{\partial x'}\right)^i
\left(\frac{\partial}{\partial y}
-\frac{\partial}{\partial y'}\right)^j
f\left(x,y\right)g\left(x',y'\right)\right|_{x'=x,\, y'=y}.
\end{equation*}
\begin{thm}\label{thm:explicit}
For a pair of solutions $\tau$, $\overline{\tau}$ to
\eqref{c_bi_1}--\eqref{c_bi_3},
we define $\theta$ and $\gamma$ by
\begin{align*}
\theta &=
2\log{\frac{\tau}{\overline{\tau}}},\\
\gamma &=
\frac{1}{2} \frac{\partial}{\partial y}
\begin{bmatrix}
- \log \left(\tau\overline{\tau}\right)\\
\log \left(\tau/\overline{\tau}\right)
\end{bmatrix}.
\end{align*}
Then,
for any $y$,
$\theta$ and $\gamma$ satisfy the potential defocusing mKdV
equation \eqref{pot_mKdV} and the defocusing mKdV flow \eqref{mkdvflow}.
\end{thm}
\noindent{\it Proof.}
First we have from \eqref{c_bi_1} and \eqref{c_bi_1-2} that
\begin{align*}
T=\gamma_x
&=\frac{1}{2}
\begin{bmatrix}
- \left(\log\tau + \log\overline{\tau}\right)_{xy}\\
\left(\log\tau - \log\overline{\tau}\right)_{xy}
\end{bmatrix}
=\frac{1}{4}
\begin{bmatrix}
-\tau^{-2}D_{x}D_{y}\tau\cdot\tau
-\overline{\tau}^{-2}D_{x}D_{y}\overline{\tau}\cdot\overline{\tau}\\
\tau^{-2}D_{x}D_{y}\tau\cdot\tau
-\overline{\tau}^{-2}D_{x}D_{y}\overline{\tau}\cdot\overline{\tau}
\end{bmatrix}\\
&=\frac{1}{2}
\begin{bmatrix}
\left(\tau/\overline{\tau}\right)^2
+ \left(\tau/\overline{\tau}\right)^{-2}\\
\left(\tau/\overline{\tau}\right)^2
- \left(\tau/\overline{\tau}\right)^{-2}
\end{bmatrix}
=\frac{1}{2}
\begin{bmatrix}
e^\theta + e^{-\theta}\\
e^\theta - e^{-\theta}
\end{bmatrix}
=\begin{bmatrix}
\cosh \theta\\
\sinh \theta
\end{bmatrix},
\end{align*}
which yields that $T$ is positive pointing,
and we have the Frenet formula \eqref{frenet}.
Thus $T_t=\theta_t N$,
by which it is sufficient for \eqref{deformation-of-frame} to show
\begin{equation}\label{c_pot}
\theta_t = \kappa_{xx}-\frac{\kappa^{3}}{2}.
\end{equation}
We have
\begin{align*}
\theta_t - \kappa_{xx} + \frac{1}{2} \kappa^3
&= 2 \frac{D_t\, \tau\cdot\overline{\tau}}{\tau\overline{\tau}}
- 2 \left(\frac{D_x\, \tau\cdot\overline{\tau}}{\tau\overline{\tau}}\right)_{xx}
+ 4 \left(\frac{D_x\, \tau\cdot\overline{\tau}}{\tau\overline{\tau}}\right)^3\\
&= \frac{2}{\tau\overline{\tau}}
\left(D_t\, \tau\cdot\overline{\tau}
    - D^3_x\, \tau\cdot\overline{\tau}
    + 3 \frac{D^2_x\, \tau\cdot\overline{\tau}}{\tau\overline{\tau}}
    D_x\, \tau\cdot\overline{\tau}\right).
\end{align*}
On the other hand,
it immediately follows from
the bilinear equations \eqref{c_bi_2} and \eqref{c_bi_3},
hence we have \eqref{c_pot}.
Equation \eqref{c_pot} is the potential defocusing mKdV
equation \eqref{pot_mKdV}.
\qqed\\
We give a solution to the bilinear equations \eqref{c_bi_1}--\eqref{c_bi_3}.
For a positive integer $N$ and an integer $k$,
we denote by $\rho_N$ the determinant
\begin{equation}\label{def:rho}
\rho_N\left(k\right)=\det
\begin{bmatrix}
{f}^{(1)}_k & {f}^{(1)}_{k+1} & \cdots & {f}^{(1)}_{k+N-1} \\
{f}^{(2)}_k & {f}^{(2)}_{k+1} & \cdots & {f}^{(2)}_{k+N-1}\\
\vdots & \vdots & \ddots & \vdots \\
{f}^{(N)}_k & {f}^{(N)}_{k+1} & \cdots & {f}^{(N)}_{k+N-1}
\end{bmatrix},
\end{equation}
and set $\rho_0 = 1$.
Here the entries $f_{k+j-1}^{(i)}$ are functions in $x$, $t$,
and auxiliary variables $y$, $z$.
\begin{prop}\label{prop:singular}
For a positive integer $N$ and an integer $k$,
define the entries of \eqref{def:rho} by
\begin{align*}
f_n^{(i)} &= \alpha_i p_i^n e^{\eta_i} + \beta_i
\left(-p_i\right)^n e^{-\eta_i},\\
\eta_i &= p_i x+4p_i^3 t+\frac{y}{p_i},
\end{align*}
where $\alpha_i, \beta_i \in \mathbb{R}$ and $p_i \in \mathbb{R}^\times$
are arbitrary constants,
and set
\begin{equation*}
\tau_N \left(k\right)
= \frac{e^{-xy}}{\prod_{i=1}^N p_i^k}\,
\rho_N\left(k\right).
\end{equation*}
Then the pair of functions
$\tau=\tau_N\left(k\right)$ and $\overline{\tau}=\tau_N\left(k+1\right)$
satisfy the bilinear equations \eqref{c_bi_1}--\eqref{c_bi_3}.
In particular,
\begin{equation}\label{theta:singular}
\theta = 2\log \frac{\rho_N\left(k\right)}{\rho_N\left(k+1\right)}
\end{equation}
satisfies the potential defocusing mKdV equation \eqref{pot_mKdV}.
\end{prop}
To prove this,
we make use of the following well-known lemma.
\begin{lem}\label{lem:dtoda}
Let $f_n^{(i)}\ \left(i, n \in \mathbb{Z}\right)$ be
a sequence of functions in $x$, $t$, $y$ and $z$, which satisfy
\begin{equation}\label{c,d_tau_f_con}
\begin{split}
\frac{\partial}{\partial x} f_n^{(i)}
&= f_{n+1}^{(i)},\quad
\frac{\partial}{\partial z} f_n^{(i)}
= f_{n+2}^{(i)},\\
\frac{\partial}{\partial t} f_n^{(i)}
&= -4 f_{n+3}^{(i)},\quad
\frac{\partial}{\partial y} f_n^{(i)}
= f_{n-1}^{(i)}.
\end{split}
\end{equation}
For a positive integer $N$ and an integer $k$,
define a function $\sigma \left(k\right)=\sigma\left(x,t;y,z;k\right)$ by
$\sigma\left(k\right)=e^{-xy} \rho_N\left(k\right)$.
Then $\sigma \left(k\right)$ satisfies the bilinear equations
\begin{align}
&D_{x}D_{y}\, \sigma\left(k\right)\cdot\sigma\left(k\right)
=-2\sigma\left(k+1\right)\sigma\left(k-1\right),\label{c_bi_1_2}\\
&\left(D_x^2-D_z\right) \sigma\left(k+1\right)\cdot\sigma\left(k\right)=0,\label{c_bi_2_2}\\
&\left(D_x^3+D_t+3D_xD_z\right) \sigma\left(k+1\right)\cdot\sigma\left(k\right)=0.\label{c_bi_3_2}
\end{align}
\end{lem}
The system of bilinear equations \eqref{c_bi_1_2}--\eqref{c_bi_3_2}
are included in the discrete two-dimensional Toda lattice hierarchy
\cite{inoguchi:2012:2,hirota:2004,jimbo:1983,maruno:1998,maruno:2006,ohta:1993,ohta:1993:2,tsujimoto:2002,ueno:1984}.
A typical example of $f_n^{(i)}$ satisfying
the condition \eqref{c,d_tau_f_con} is given by
\begin{equation}\label{c,d_tau_f_1}
f_n^{(i)} =
\alpha_i p_i^n e^{p_i x + p_i^2 z -4 p_i^3 t + p_i^{-1} y}+
\beta_i q_i^n e^{q_i x + q_i^2 z -4 q_i^3 t + q_i^{-1} y},
\end{equation}
where $\alpha_i, \beta_i \in \mathbb{R}$
and $p_i, q_i \in \mathbb{R}^\times$ are arbitrary constants.
\\\noindent{\it Proof.}
Let us prove Proposition \ref{prop:singular}.
Fix $N$ and $k$.
Let $\sigma\left(k\right)$ a function as introduced in Lemma \ref{lem:dtoda}
with entries \eqref{c,d_tau_f_1}.
Imposing on it the reduction condition $q_i=-p_i$,
we have for all integers $i$ and $n$ that
\begin{equation*}
f_{n+2}^{(i)}
= \frac{\partial}{\partial z} f_n^{(i)}
= {p_i}^2 f_n^{(i)},
\end{equation*}
which yields
\begin{equation}\label{d2dt_impose_3_4}
\sigma\left(k+2\right)=\sigma\left(k\right)\prod_{i=1}^N p_i^2,\quad
\frac{\partial}{\partial z}\sigma\left(k\right)=
\sigma\left(k\right) \sum_{i=1}^N p_i^2.
\end{equation}
Because $\sigma\left(k\right)=\tau_N\left(k\right)
\prod_{i=1}^N p_i^k$,
equations \eqref{d2dt_impose_3_4} are rewritten
in terms of $\tau_N\left(k\right)$ as
\begin{equation}\label{d2dt_impose_4_1}
\tau_N\left(k+2\right)=\tau_N\left(k\right),\quad
\frac{\partial}{\partial z}\tau_N\left(k\right)
=\tau_N\left(k\right)\sum_{i=1}^N p_i^2.
\end{equation}
Therefore the system of bilinear equations \eqref{c_bi_1}--\eqref{c_bi_3} with $c=0$
immediately follows from \eqref{c_bi_1_2}--\eqref{c_bi_3_2}
on writing $\tau=\tau_N\left(k\right)$ and
$\overline{\tau}=\tau_N\left(k+1\right)$.
\qqed\\
\begin{figure}[H]
 \centering
 \subfigure[]{
 \includegraphics[width=0.3\linewidth]{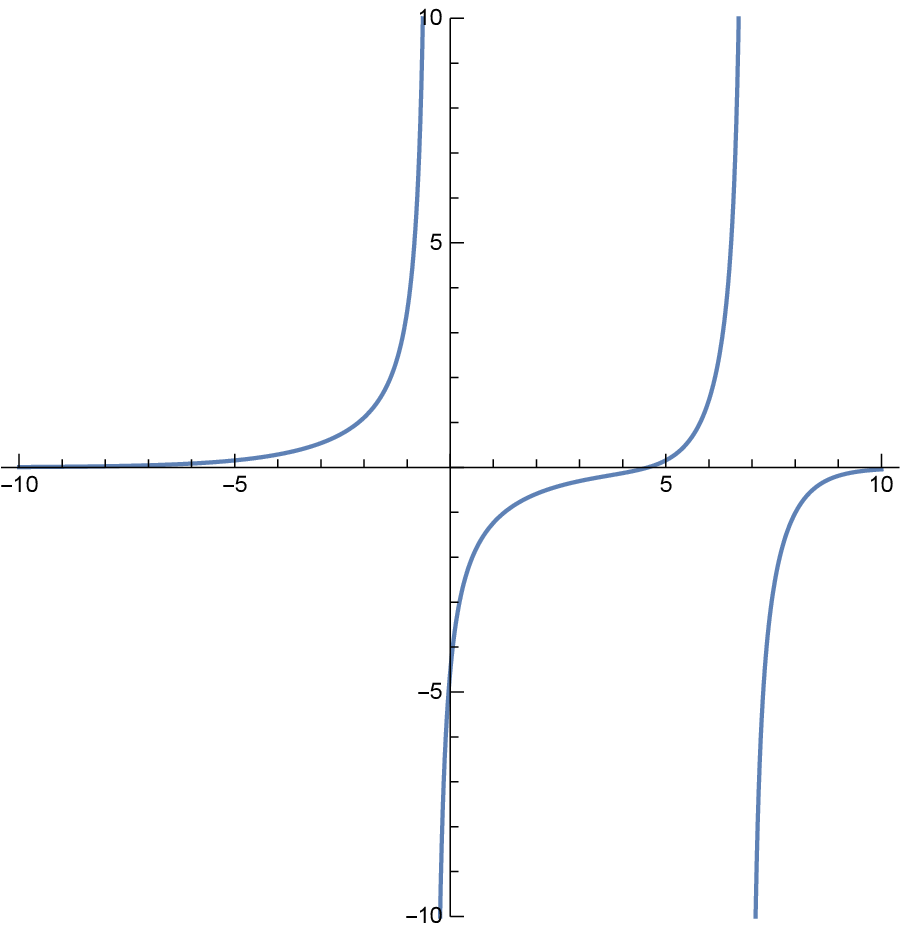}
 }
 \centering
 \subfigure[]{
 \includegraphics[width=0.3\linewidth]{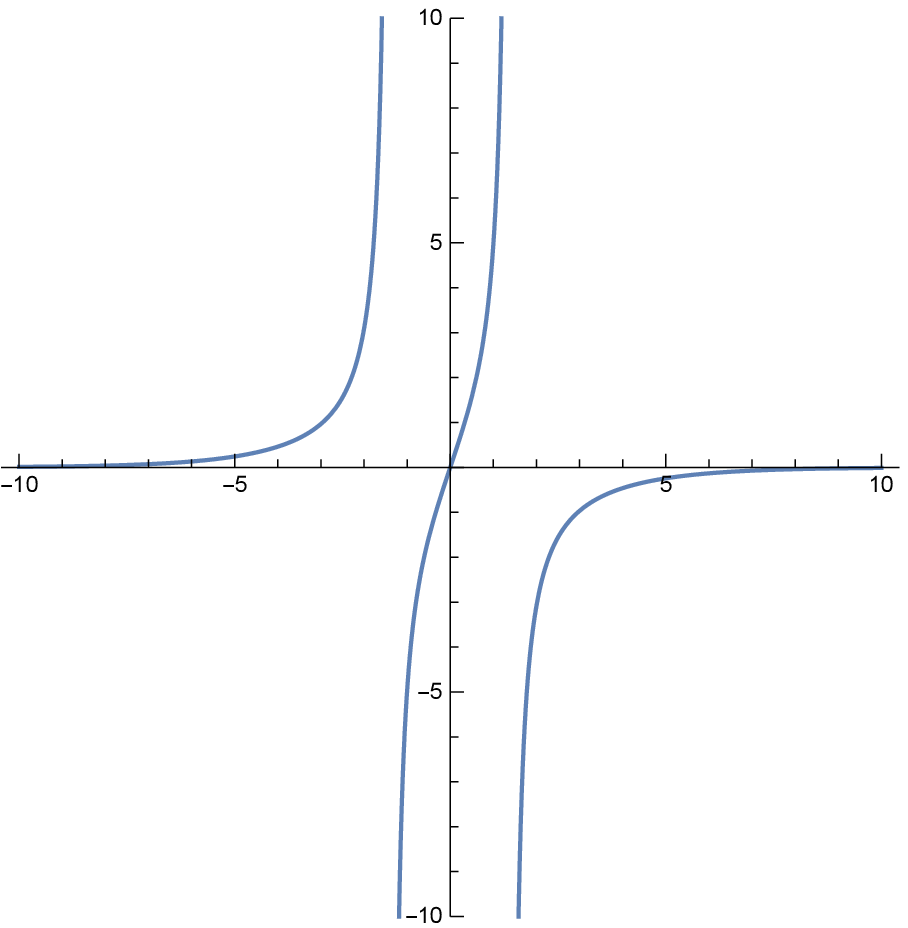}
 }
 \centering
 \subfigure[]{
 \includegraphics[width=0.3\linewidth]{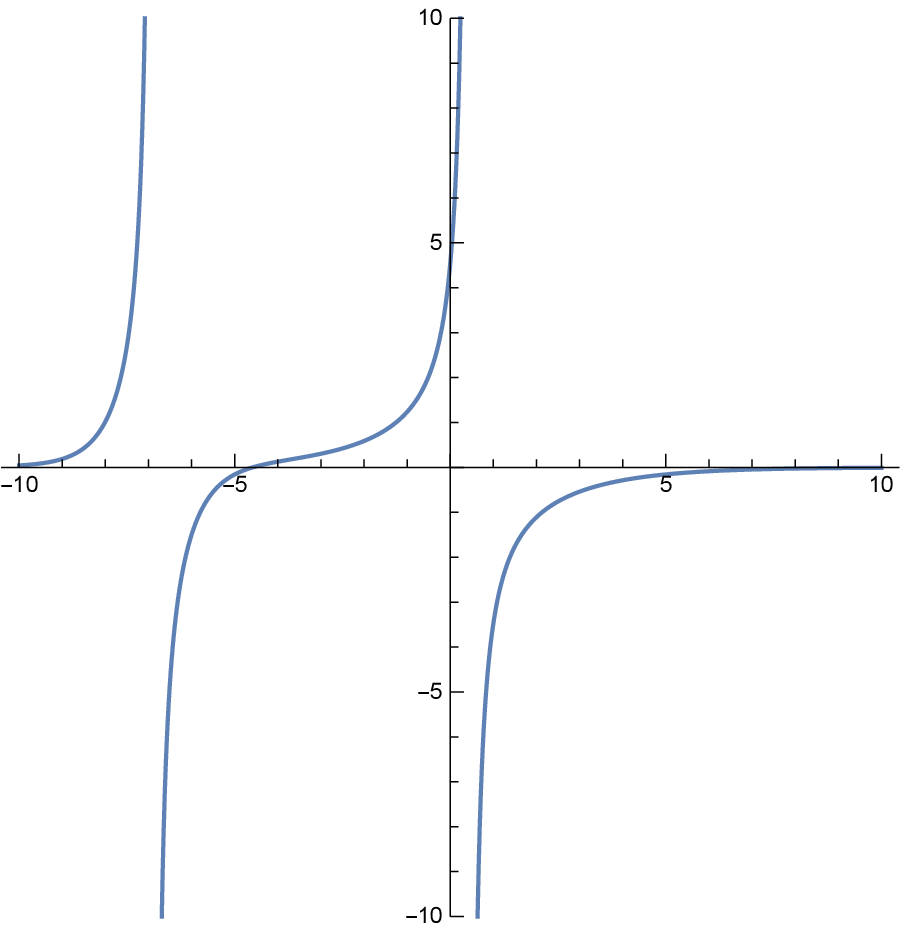}
 }
 \caption{Profiles of the solution to the defocusing mKdV equation \eqref{mkdv}
 $\kappa=\theta_x$, where $\theta$ is in \eqref{theta:singular}.
 Parameters in Proposition \ref{prop:singular} are
 $N=2$, $p_1=0.3$, $p_2=0.9$, $\alpha_1=\beta_1=\beta_2=1$ and $\alpha_2=-1$, and
 $t=-7$ (a), $t=0$ (b), $t=7$ (c).}
 \label{figure:mKdV}
\end{figure}
\begin{figure}[H]
 \centering
 \subfigure[]{
 \includegraphics[width=0.3\linewidth]{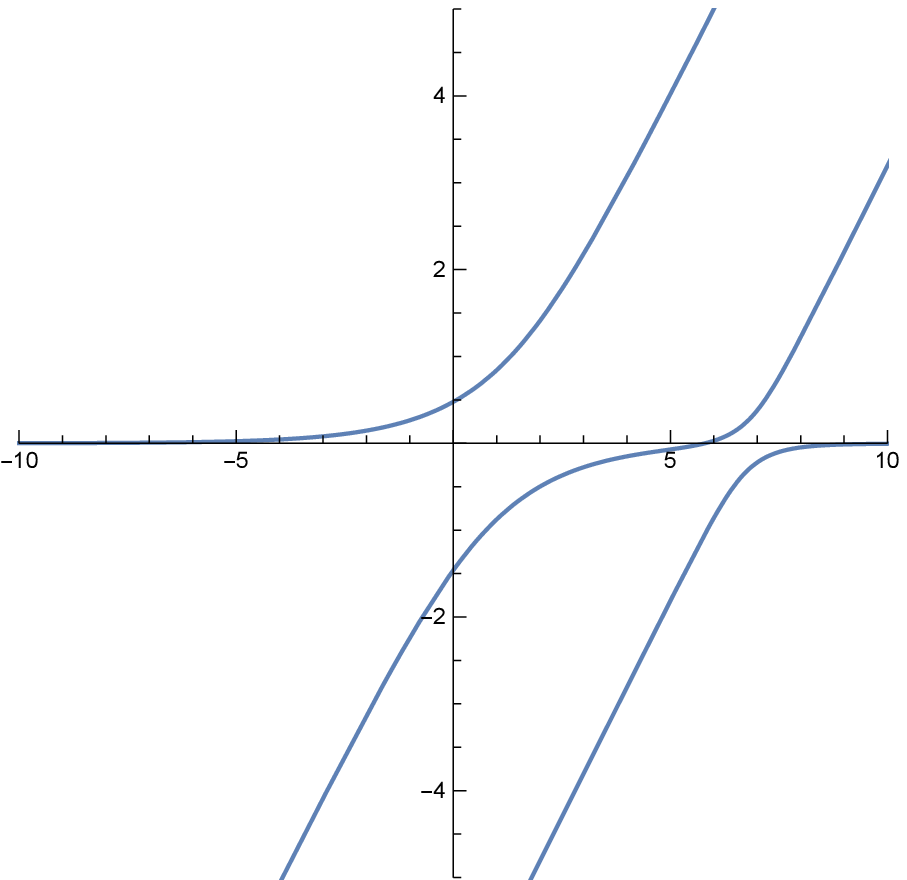}
 }
 \centering
 \subfigure[]{
 \includegraphics[width=0.3\linewidth]{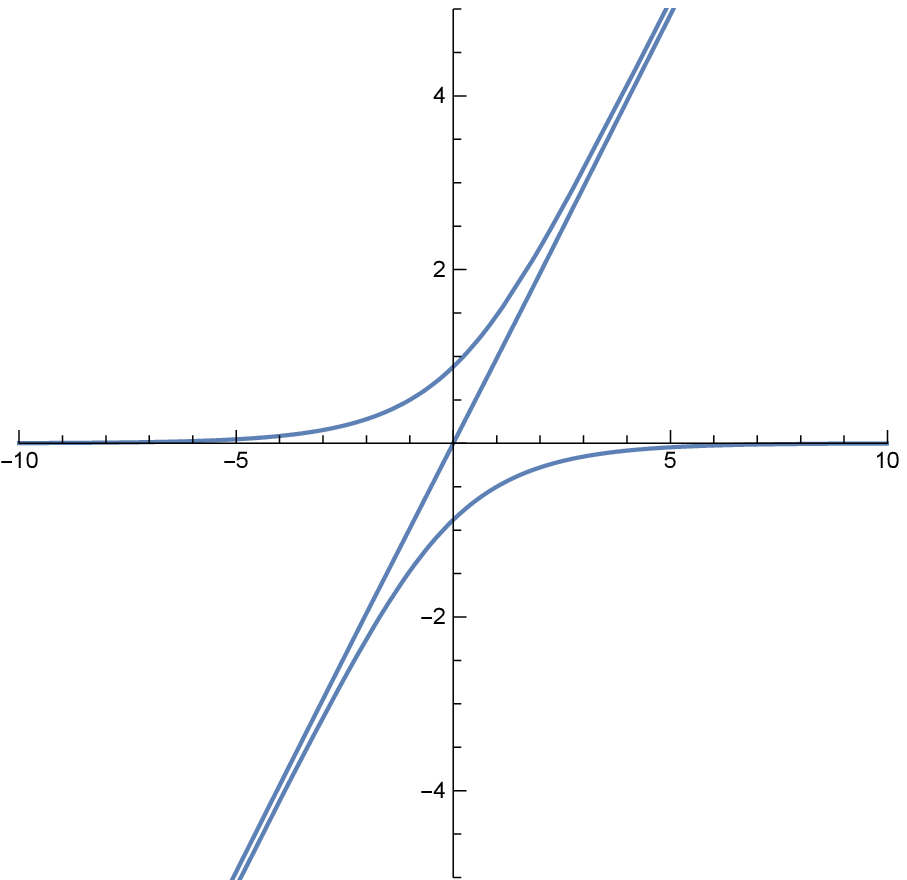}
 }
 \centering
 \subfigure[]{
 \includegraphics[width=0.3\linewidth]{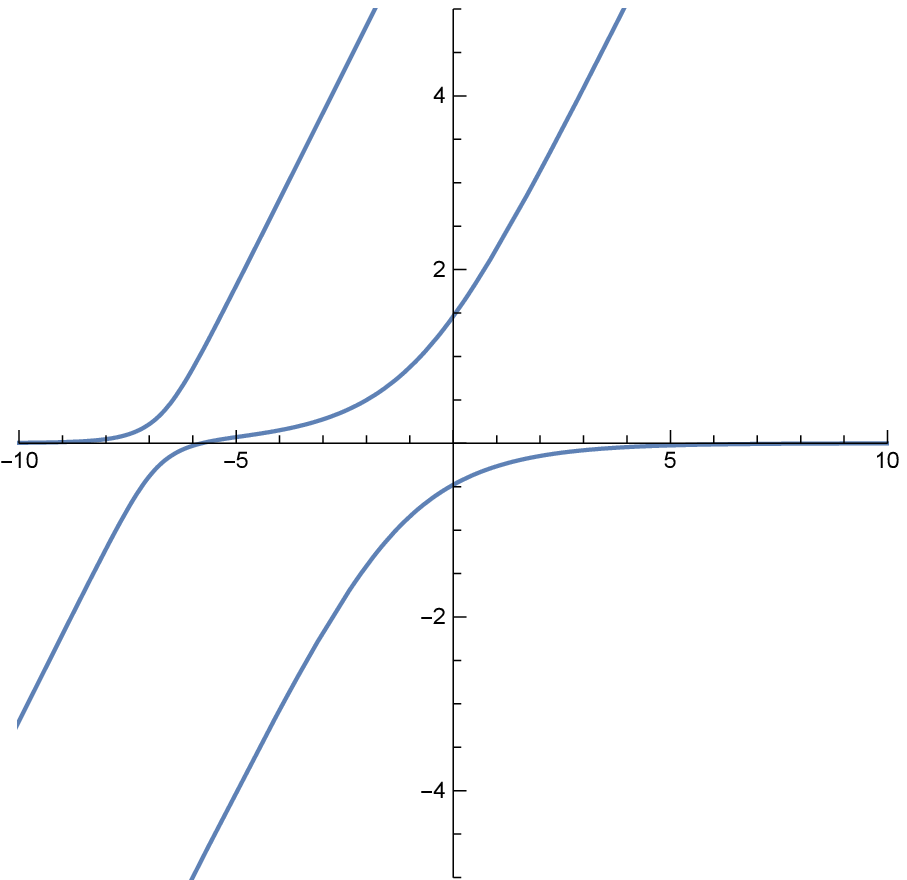}
 }
 \caption{Profiles of the defocusing mKdV flow \eqref{mkdvflow}
 given in Theorem \ref{thm:explicit}.
 Parameters in Proposition \ref{prop:singular} are
 $N=2$, $p_1=0.3$, $p_2=0.9$, $\alpha_1=\beta_1=\beta_2=1$ and $\alpha_2=-1$, and
 $t=-7$ (a), $t=0$ (b), $t=7$ (c).}
 \label{figure:singular}
\end{figure}
\subsection{Regular solutions}\label{section:solution2}

The formula in Theorem \ref{thm:explicit}
together with Proposition \ref{prop:singular}
gives exact solutions to
the potential defocusing mKdV equation \eqref{pot_mKdV}
and the corresponding defocusing mKdV flow \eqref{mkdvflow}.
However,
as shown in Figures \ref{figure:mKdV} and Figure \ref{figure:singular},
they have singular points
since the $\tau$ function $\rho_N \left(k\right)$, $\rho_N \left(k+1\right)$
have zeros in general.
Actually, it seems difficult to choose parameters such that
those $\tau$ functions are positive valued.
In contrast to this,
by using the $\tau$ functions
$\rho_N \left(k\right)$, $\rho_{N+1} \left(k\right)$ appropriately,
we are able to give regular solutions
to the potential defocusing mKdV equation \eqref{pot_mKdV}
according to the idea suggested in
\cite{ohta:1989}.
With this choice of $\tau$ functions, however,
it is difficult to construct an explicit formula for the
corresponding motion of curves on the Minkowski plane
since functions $\tau=\tau_N(k)$ and $\bar{\tau}=\tau_{N+1}(k)$,
which are defined by the same way as that in Proposition \ref{prop:singular},
do not satisfy the bilinear equations \eqref{c_bi_1} and \eqref{c_bi_1-2}.
Therefore, we only present a regular solutions to
the potential defocusing mKdV equation
and observe the corresponding defocusing mKdV flow by using a numerical method.
First we introduce the following lemma:
\begin{lem}\label{lem:mkp}
Let $g_n^{(i)}\ \left(i, n \in \mathbb{Z}\right)$ be
a sequence of functions in ${u}$, ${v}$ and ${w}$, which satisfy
\begin{equation}\label{mKP_f}
\frac{\partial}{\partial {u}} g_n^{(i)} =g_{n+1}^{(i)},\quad
\frac{\partial}{\partial {v}} g_n^{(i)} =g_{n+2}^{(i)},\quad
\frac{\partial}{\partial {w}} g_n^{(i)} =g_{n+3}^{(i)}.
\end{equation}
For a positive integer $N$ and an integer $k$,
we denote by $h_N$ the determinant
\begin{equation*}
h_N\left(k\right)=\det
\begin{bmatrix}
g^{(1)}_k & g^{(1)}_{k+1} & \cdots & g^{(1)}_{k+N-1} \\
g^{(2)}_k & g^{(2)}_{k+1} & \cdots & g^{(2)}_{k+N-1}\\
\vdots & \vdots & \ddots & \vdots \\
g^{(N)}_k & g^{(N)}_{k+1} & \cdots & g^{(N)}_{k+N-1}
\end{bmatrix},
\end{equation*}
and set $h_0 = 1$.
Then $h_N$ satisfies the bilinear equations of the mKP hierarchy
\cite{jimbo:1983,freeman:1983,gesztesy:1991,hirota:1988,hirota:1988:2}
\begin{equation}\label{mkphierarchy}
\begin{split}
& \left(D_{u}^2 - D_{v}\right) h_{N+1}(k) \cdot h_{N}(k)=0,\\
& \left(D_{u}^3 -4 D_{w} +3 D_{u} D_{v}\right) h_{N+1}(k) \cdot h_{N}(k)=0.
\end{split}
\end{equation}
\end{lem}

Here we note that,
for integers $n$, $i$ and parameters $\alpha_i, \beta_i, p_i, q_i$,
the function
\begin{equation*}
g_n^{(i)} \left({u},{v},{w}\right)
= \alpha_i p_i^n e^{p_i {u} + p_i^2 {v} + p_i^3 {w}}
+ \beta_i q_i^n e^{q_i {u} + q_{i}^2 {v} + q_i^3 {w}}
\end{equation*}
satisfies \eqref{mKP_f}.
Imposing on this a reduction condition $q_i=-p_i$,
it follows that the determinant function $h_N$ satisfies
$\left(\partial/\partial {v}\right) h_N
= h_N \sum_{i=1}^N p_i^2$,
and hence $D_{v}\, h_{N+1}\cdot h_N=\left(p_{N+1}\right)^2 h_{N+1}h_N$.
Thus the bilinear equations \eqref{mkphierarchy} become
\begin{equation}\label{mkphierarchy-reduction}
\begin{split}
& \left(D_{u}^2 - {p_{N+1}}^2 \right) h_{N+1}\cdot h_{N}=0,\\
& \left(D_{u}^3 -4 D_{w} + 3 {p_{N+1}}^2 D_{u}\right) h_{N+1}\cdot h_{N}=0.
\end{split}
\end{equation}
Using Lemma \ref{lem:mkp} and
the bilinear equations \eqref{mkphierarchy-reduction},
we construct a class of regular solutions
the potential defocusing mKdV equation \eqref{pot_mKdV}
by applying the Galilean transformation
and by imposing appropriate conditions on the parameters as follows.
We denote by $\mathrm{sgn}$ the sign function:
\begin{equation*}
\sgn{x} =
\begin{cases}
1 & \text{if $x>0$}\\
-1 & \text{if $x<0$}.
\end{cases}
\end{equation*}
\begin{prop}\label{prop:regular}
Fix an integer $k$.
For a positive integer $N$,
define the entries of \eqref{def:rho} by
\begin{align*}
f_n^{(i)} & = \alpha_i p_i^n e^{\eta_i}
+\beta_i \left(-p_i\right)^n e^{-\eta_i},\\
\eta_i & = p_i x
+4{p_i}^3 \left(1-\frac{3}{2} a {p_i}^2\right)t,
\end{align*}
where $a$, $\alpha_i$, $\beta_i$ and $p_i$ are arbitrary real constants
for $i=1,2,\ldots,N+1$.
Then putting $a=p_{N+1}^{-2}$, the function
\begin{equation}\label{theta:regular}
\theta = 2\log \frac{\rho_N\left(k\right)}{\rho_{N+1}\left(k\right)}
\end{equation}
satisfies the potential defocusing mKdV equation \eqref{pot_mKdV}.
Moreover,
if we choose the parameters
in such a way that
\begin{equation}\label{parameter:regular}
  \begin{split}
    0 < p_1 < \cdots < p_{N+1}, \\
    \alpha_i > 0,\quad
    \sgn{\beta_i} = \left(-1\right)^{k+i-1}
  \end{split}
\end{equation}
for all $i$,
then $\theta$ gives a regular solution to \eqref{pot_mKdV}.
\end{prop}
\noindent{\it Proof.}\label{prop:regular:sol}
First we only prove that $\theta$ is a solution to
the potential defocusing mKdV equation \eqref{mkdv},
and afterward we shall verify regularity of $\theta$
on the condition \eqref{parameter:regular}.
We change the independent variables
from $\left({u},{v},{w}\right)$ to $\left(x,t,z\right)$ by
\begin{equation*}
x = u + \frac{3}{2} r^2 w,\quad
t = \frac{1}{4} w,\quad
z = {v},
\end{equation*}
where $r$ is a real constant,
and define
$\rho_N\left(x,t\right) = h_N\left({u},0,{w}\right)$.
Then putting $r = p_{N+1}^2$,
we have
\begin{align*}
& \left(D_x^2- r^2 \right) \rho_{N+1}\cdot\rho_N=0,\\
& \left(D_x^3-D_t-3 r^2 D_x\right) \rho_{N+1}\cdot\rho_N=0.
\end{align*}
Thus the pair $\rho_{N+1},\,\rho_N$ gives a solution to the bilinear equations
\eqref{c_bi_2} and \eqref{c_bi_3}.
Therefore \eqref{theta:regular} satisfies
the potential defocusing mKdV equation \eqref{pot_mKdV}.
\qqed\\
\begin{figure}[H]
 \centering
 \subfigure[]{
 \includegraphics[width=0.3\linewidth]{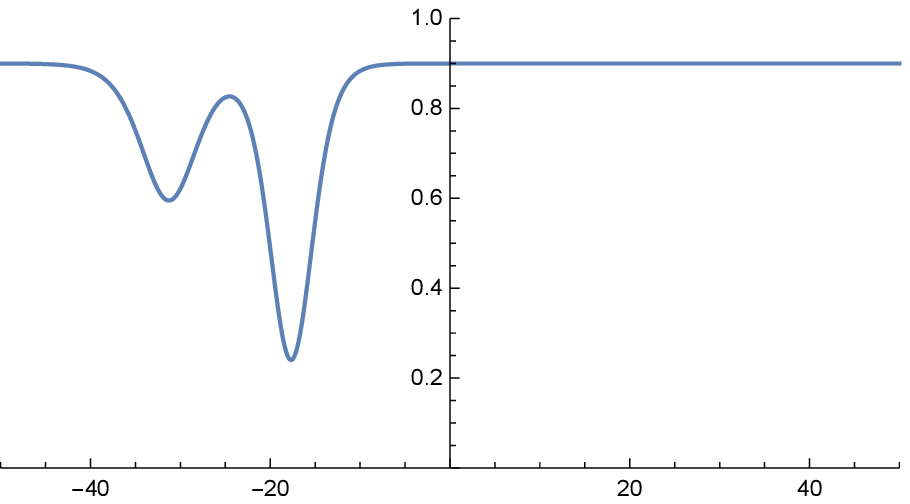}
 }
 \centering
 \subfigure[]{
 \includegraphics[width=0.3\linewidth]{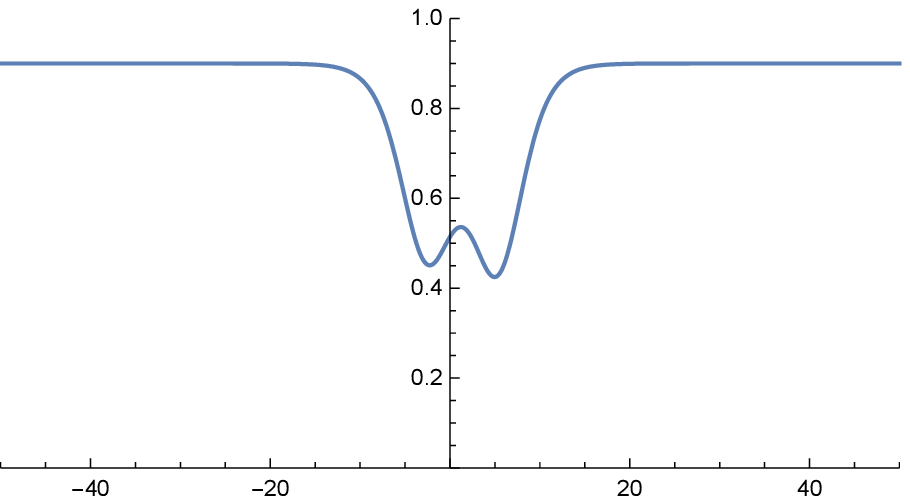}
 }
 \centering
 \subfigure[]{
 \includegraphics[width=0.3\linewidth]{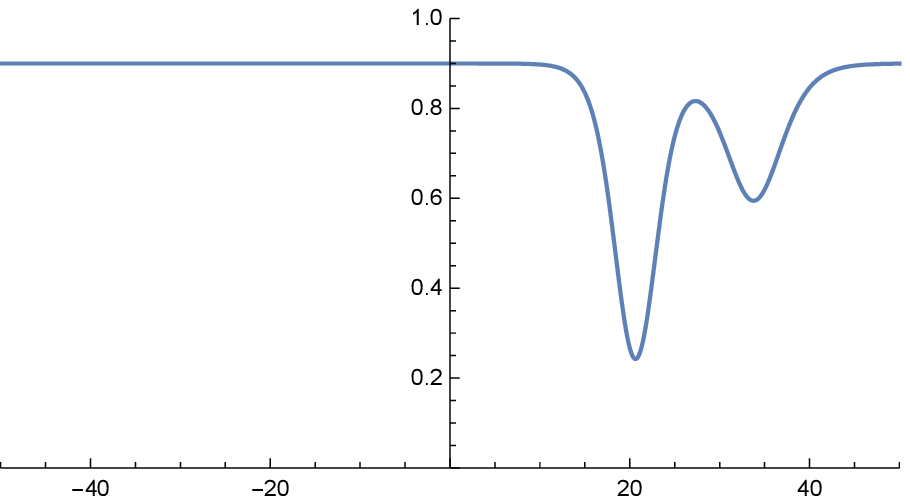}
 }
 \caption{Profiles of the solution to the defocusing mKdV equation \eqref{mkdv}
 $\kappa=\theta_x$, where $\theta$ is in \eqref{theta:regular}.
 Parameters in Proposition \ref{prop:regular} are
 $N=2$, $p_1=0.5$, $p_2=0.7$, $p_3=0.9$, $\alpha_1=\alpha_2=\alpha_3=\beta_1=1$,
 $\beta_2=-1$ and $\beta_3=0$, and
 $t=-30$ (a), $t=0$ (b), $t=30$ (c).}
 \label{figure:dark}
\end{figure}
\begin{figure}[H]
 \centering
 \subfigure[]{
 \includegraphics[width=0.3\linewidth]{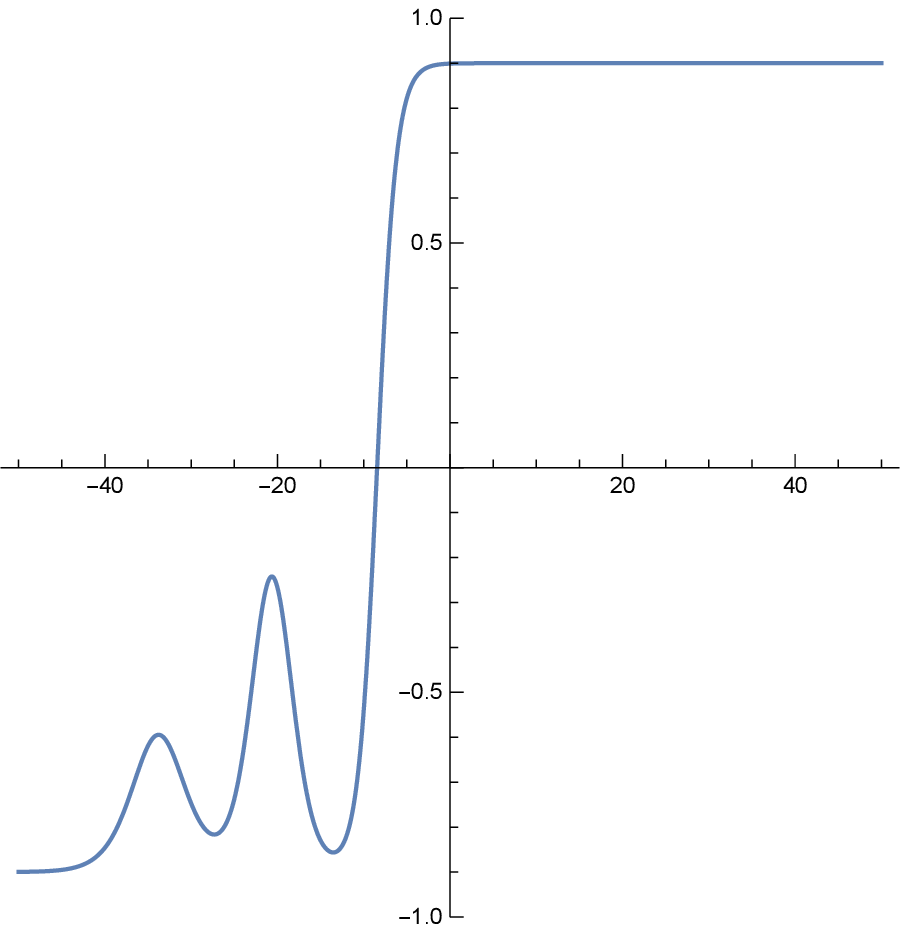}
 }
 \centering
 \subfigure[]{
 \includegraphics[width=0.3\linewidth]{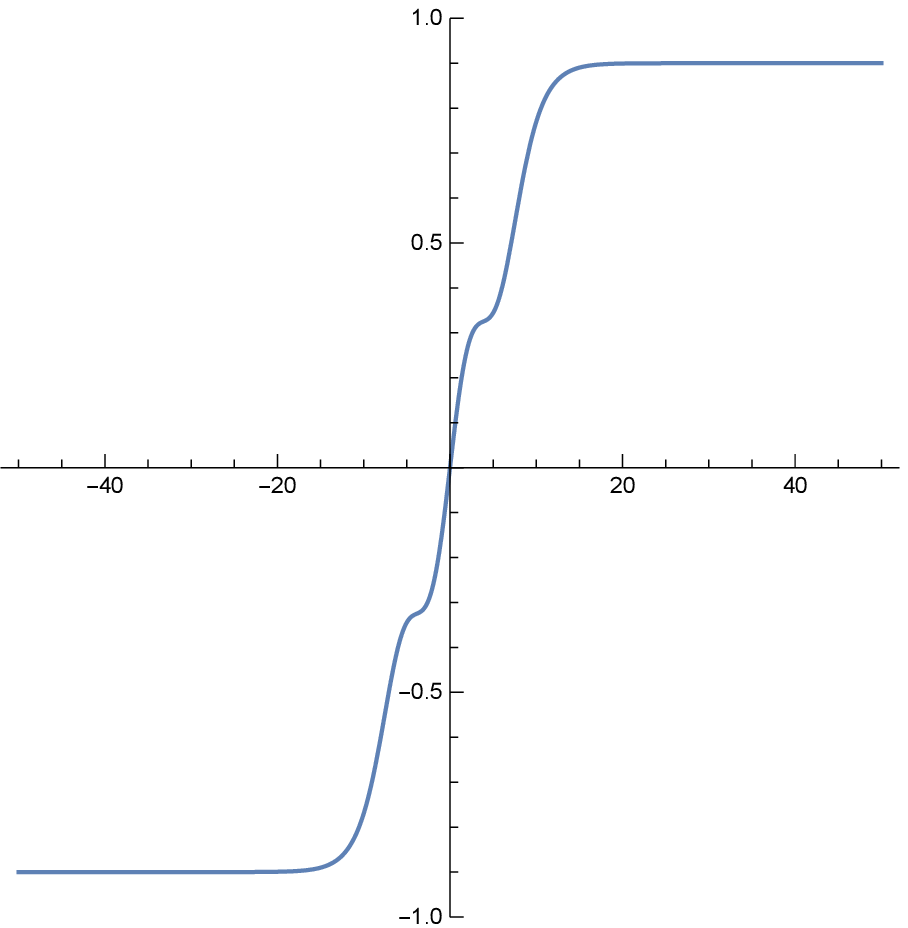}
 }
 \centering
 \subfigure[]{
 \includegraphics[width=0.3\linewidth]{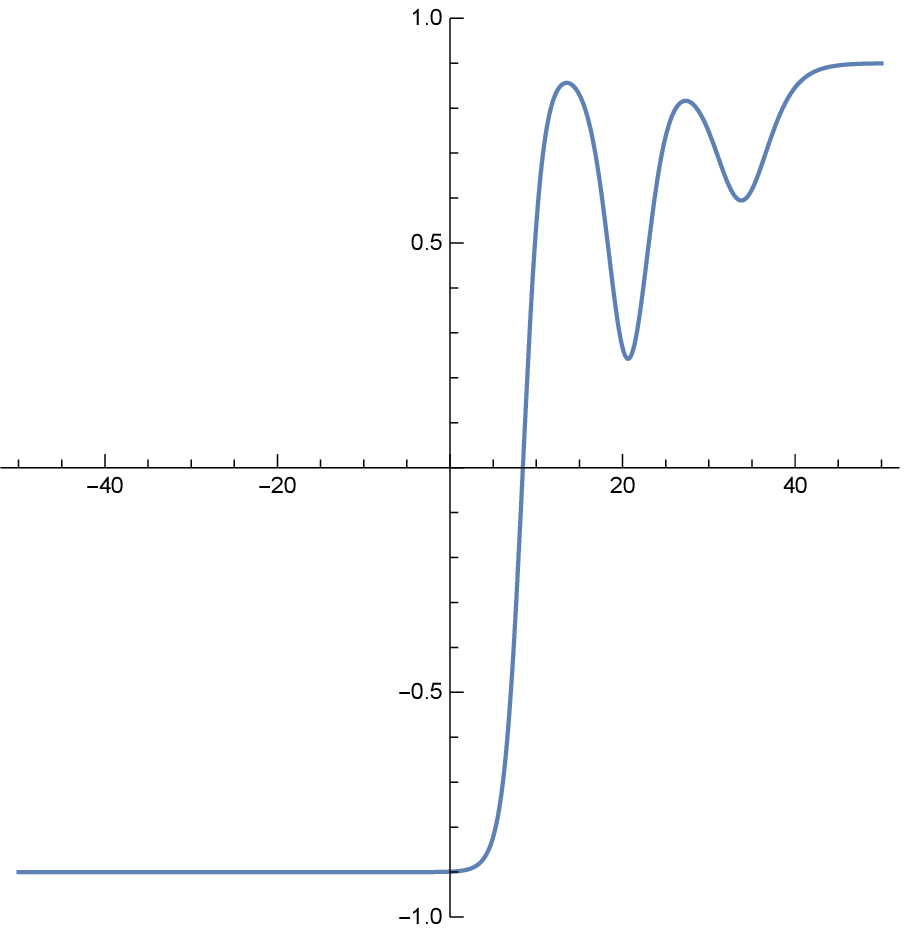}
 }
 \caption{Profiles of the solution to the defocusing mKdV equation \eqref{mkdv}
 $\kappa=\theta_x$, where $\theta$ is in \eqref{theta:regular}.
 Parameters in Proposition \ref{prop:regular} are
 $N=2$, $p_1=0.5$, $p_2=0.7$, $p_3=0.9$, $\alpha_1=\alpha_2=\alpha_3=\beta_1=\beta_3=1$
 and $\beta_2=-1$, and
 $t=-30$ (a), $t=0$ (b), $t=30$ (c).}
 \label{figure:kink}
\end{figure}
Now we establish the regularity of $\theta$ by proving the positivity of $\rho_N$
in Proposition \ref{prop:regular:sol}.
The condition on the parameters \eqref{parameter:regular} plays a crucial role.
\begin{lem}\label{lem:positivity}
Fix an integer $k$.
For a positive integer $N$,
let $F$ be an $N \times N$ matrix
\begin{align*}
F &=
\begin{bmatrix}
  \alpha_{1} p_{1}^{k} e^{\eta_{1}}
  + \beta_{1} q_{1}^{k} e^{\xi_{1}} & \cdots &
  \alpha_{1} p_{1}^{k+N-1} e^{\eta_{1}}
  + \beta_{1} q_{1}^{k+N-1} e^{\xi_{1}}\\
  \vdots & \ddots & \vdots \\
  \alpha_{N} p_{N}^{k} e^{\eta_{N}}
  + \beta_{N} q_{N}^{k} e^{\xi_{N}} & \cdots &
  \alpha_{N} p_{N}^{k+N-1} e^{\eta_{N}}
  + \beta_{N} q_{N}^{k+N-1} e^{\xi_{N}}
\end{bmatrix},
\end{align*}
where $\alpha_i$, $\beta_i$, $p_i$, $q_i$, $\eta_i$ and
$\xi_i$ are arbitrary real parameters for all $i$.
Choose these parameters as
{
\begin{equation}\label{eqn:param_regular_tau}
\begin{split}
q_N < q_{N-1} < \cdots < q_1 < 0 < p_1 < p_2 < \cdots < p_N, \\
\alpha_i > 0,\quad \sgn{\beta_i} = {(-1)}^{k+i-1}\quad \left( i=1,2,\ldots,N \right).
\end{split}
\end{equation}
}
Then $\det F$ is positive.
\end{lem}
\noindent{\it Proof.}
The matrix $F$ is expressed as a product of
$N \times 2N$ and $2N \times N$ matrices as
\begin{equation*}
F =
  \begin{bmatrix}
    \alpha_1 & \beta_1 & 0 & 0 & \cdots & 0 & 0 \\
    0 & 0 & \alpha_2 & \beta_2 & \cdots & 0 & 0 \\
    \vdots & \vdots & \vdots & \vdots & \ddots & \vdots & \vdots \\
    0 & 0 & 0 & 0 & \cdots & \alpha_N & \beta_N
  \end{bmatrix}
  \begin{bmatrix}
    p_1^k e^{\eta_1} & \cdots & p_1^{k+N-1} e^{\eta_1} \\
    q_1^k e^{\xi_1} & \cdots & q_1^{k+N-1} e^{\xi_1} \\
    p_2^k e^{\eta_2} & \cdots & p_2^{k+N-1} e^{\eta_2} \\
    q_2^k e^{\xi_2} & \cdots & q_2^{k+N-1} e^{\xi_2} \\
    \vdots & \ddots & \vdots \\
    p_N^k e^{\eta_N} & \cdots & p_n^{k+N-1} e^{\eta_N} \\
    q_N^k e^{\xi_N} & \cdots & q_n^{k+N-1} e^{\xi_N} \\
  \end{bmatrix}.
\end{equation*}
Then by using the Cauchy-Binet formula, we have
\begin{equation}\label{Cauchy-Binet}
\det F = \sum_{\mu_1, \mu_2, \ldots, \mu_N}
d \left(\mu_1, \mu_2, \ldots, \mu_N\right),
\end{equation}
{where the summation in \eqref{Cauchy-Binet} is taken over all possible
combinations of
\begin{equation*}
\mu_1 \in \left\{\alpha_1, \beta_1\right\},\ %
\mu_2 \in \left\{\alpha_2, \beta_2\right\},\ %
\ldots,\ %
\mu_N \in \left\{\alpha_N, \beta_N\right\},
\end{equation*}
}
\begin{align*}
d \left(\mu_1, \mu_2, \ldots, \mu_N\right)
&= \det
\begin{bmatrix}
\mu_1 & 0 & \cdots & 0\\
0 & \mu_2 & \cdots & 0\\
\vdots & \vdots & \ddots & \vdots\\
0 & 0 & \cdots & \mu_N
\end{bmatrix}
\det
\begin{bmatrix}
\nu_1^k e^{\omega_1} & \cdots & \nu_1^{k+N-1} e^{\omega_1} \\
\nu_2^k e^{\omega_2} & \cdots & \nu_2^{k+N-1} e^{\omega_2} \\
\vdots & \ddots & \vdots \\
\nu_N^k e^{\omega_N} & \cdots & \nu_N^{k+N-1} e^{\omega_N} \\
\end{bmatrix}\\
&=
\prod_{n=1}^N \mu_n \nu_n^k e^{\omega_n}
\det
\begin{bmatrix}
1 &  \nu_1  & \cdots & \nu_1^{N-1} \\
1 &  \nu_2  & \cdots & \nu_2^{N-1} \\
\vdots & \vdots & \ddots & \vdots \\
1 &  \nu_N  & \cdots & \nu_N^{N-1}
\end{bmatrix}\\
&=
\prod_{n=1}^N \mu_n \nu_n^k e^{\omega_n}
\prod_{1 \leq i < j \leq N} \left(\nu_j - \nu_i\right),
\end{align*}
and each $\left(\nu_i, \omega_i\right)$ is accordingly given by
\begin{equation*}
\left(\nu_i, \omega_i\right) =
\begin{cases}
\left(p_i, \eta_i\right) & \text{if $\mu_i = \alpha_i$}\\
\left(q_i, \xi_i\right) & \text{if $\mu_i = \beta_i$}.
\end{cases}
\end{equation*}
It suffices for the positivity of $\det F$ to show that
$d \left(\mu_1, \mu_2, \ldots, \mu_N\right) > 0$
for each $\left(\mu_1, \mu_2, \ldots, \mu_N\right)$.
We arbitrarily fix a choice $\left(\mu_1, \mu_2, \ldots, \mu_N\right)$.
If $\mu_i = \alpha_i$ for all $i$,
then $d \left(\mu_1, \mu_2, \ldots, \mu_N\right)$ is obviously positive.
If $\mu_i = \beta_i$ for all $i$,
then $d \left(\mu_1, \mu_2, \ldots, \mu_N\right)$ is positive,
because
\begin{equation*}
\sgn{\prod_{n=1}^N \beta_n q_n^k}
= \left(-1\right)^{N \left(N-1\right)/2},\quad
\sgn{\prod_{1\leq i<j\leq N} \left(q_j - q_i\right)}
= \left(-1\right)^{N \left(N-1\right)/2}.
\end{equation*}
For the remaining cases,
we divide the set
$\mu = \left\{\mu_1, \mu_2, \ldots, \mu_N\right\}$
into $\mu = \mu^{\alpha} \cup \mu^{\beta}$, where
\begin{align*}
\mu^\alpha =
\left\{\alpha_{j_1}, \alpha_{j_2}, \ldots, \alpha_{j_{N-r}}\right\},\quad
\mu^\beta =
\left\{\beta_{i_1}, \beta_{i_2}, \ldots, \beta_{i_r}\right\}
\end{align*}
with some $r \in \left\{1, 2, \ldots, N-1\right\}$.
Here the indices are sorted in ascending order,
namely $j_1 < j_2 < \cdots < j_{N-r}$
and $i_1 < i_2 < \cdots < i_r$.
Then we have
\begin{equation}\label{sign1}
\sgn{\prod_{n=1}^{N} \mu_n \nu_n^{k}}
= {(-1)}^{i_1 + \cdots + i_r - r},
\end{equation}
because
\begin{align*}
\prod_{n=1}^{N} \mu_n \nu_n^{k}
&= \mu_1 \mu_2 \cdots \mu_N
\left(\nu_1 \nu_2 \cdots \nu_N\right)^{k}\\
&= \alpha_{j_1} \cdots \alpha_{j_{N-r}}
\beta_{i_1} \cdots \beta_{i_r}
\left( p_{j_1} \cdots p_{j_{N-r}} q_{i_1} \cdots q_{i_r} \right)^{k}\\
&= \alpha_{j_1} \cdots \alpha_{j_{N-r}}
{\left(-1\right)}^{i_1 + \cdots + i_r + r k - r}
\left|\beta_{i_1} \cdots \beta_{i_r}\right|
{\left(-1\right)}^{rk}
\left( p_{j_1} \cdots p_{j_{n-r}}
\left|q_{i_1} \cdots q_{i_r}\right|\right)^{k}\\
&= {\left(-1\right)}^{i_1 + \cdots + i_r - r}
\prod_{n=1}^{N} \left|\mu_n \nu_n^{k}\right|.
\end{align*}
Similarly we have
\begin{equation}\label{sign2}
\sgn{\prod_{1 \leq i < j \leq N} \left(\nu_j - \nu_i\right)} =
{\left(-1\right)}^{i_1 + \cdots + i_r - r}.
\end{equation}
In fact,
for each $j \in \left\{2, 3, \ldots, N \right\}$,
it follows for all $i = 1, 2, \ldots, j-1$ that
\begin{align*}
\sgn{\left(\nu_j - \nu_i\right)}
&=
\begin{cases}
1 & \text{if $\nu_j = p_j$}\\
-1 & \text{if $\nu_j = q_j$},
\end{cases}
\end{align*}
and consequently
\begin{align*}
\sgn{\prod_{i=1}^{j-1} \left(\nu_j - \nu_i\right)}
&=
\begin{cases}
1 & \text{if $j \in \left\{j_1, \ldots, j_{N-r}\right\}$}\\
\left(-1\right)^{j-1} &
\text{if $j \in \left\{i_1, \ldots, i_r\right\}$}.
\end{cases}
\end{align*}
Thus we readily have \eqref{sign2} if $i_1 \neq 1$.
We note that,
if $i_1 = 1$,
we have
\begin{align*}
\sgn{\prod_{1 \leq i < j \leq N} \left(\nu_j - \nu_i\right)}
= \left(-1\right)^{i_2 + i_3 + \cdots + i_r - \left(r-1\right)}
= \left(-1\right)^{i_ 1 + i_2 + \cdots + i_r - r}.
\end{align*}
Therefore it follows from \eqref{sign1} and \eqref{sign2} that
\begin{align*}
\sgn{d \left(\mu_1, \mu_2, \ldots, \mu_N\right)}
=
\sgn{\left( \prod_{n=1}^N \mu_n \nu_n^k e^{\omega_n}
\prod_{1 \leq i < j \leq N} \left(\nu_j - \nu_i\right) \right)}
=
1.
\end{align*}
Thus every $d \left(\mu_1, \mu_2, \ldots, \mu_N\right)$ is positive,
and hence $\det F$ is positive.
\qqed\\
The positivity of $\rho_N(k)$ and thus the regularity of $\theta$ follow
immediately from Lemma \ref{lem:positivity} by putting
$q_i = -p_i$ and $\xi_i = -\eta_i$ for $i=1,2,\ldots,N$.
We illustrate some spacelike curves
by using the representation formula \eqref{representation}
and the regular solution in Proposition \ref{prop:regular}.
By applying a numerical integration in \eqref{representation}
with $\theta$ in \eqref{theta:regular},
we have the following figures.
\begin{figure}[H]
 \centering
 \subfigure[]{
 \includegraphics[width=0.3\linewidth]{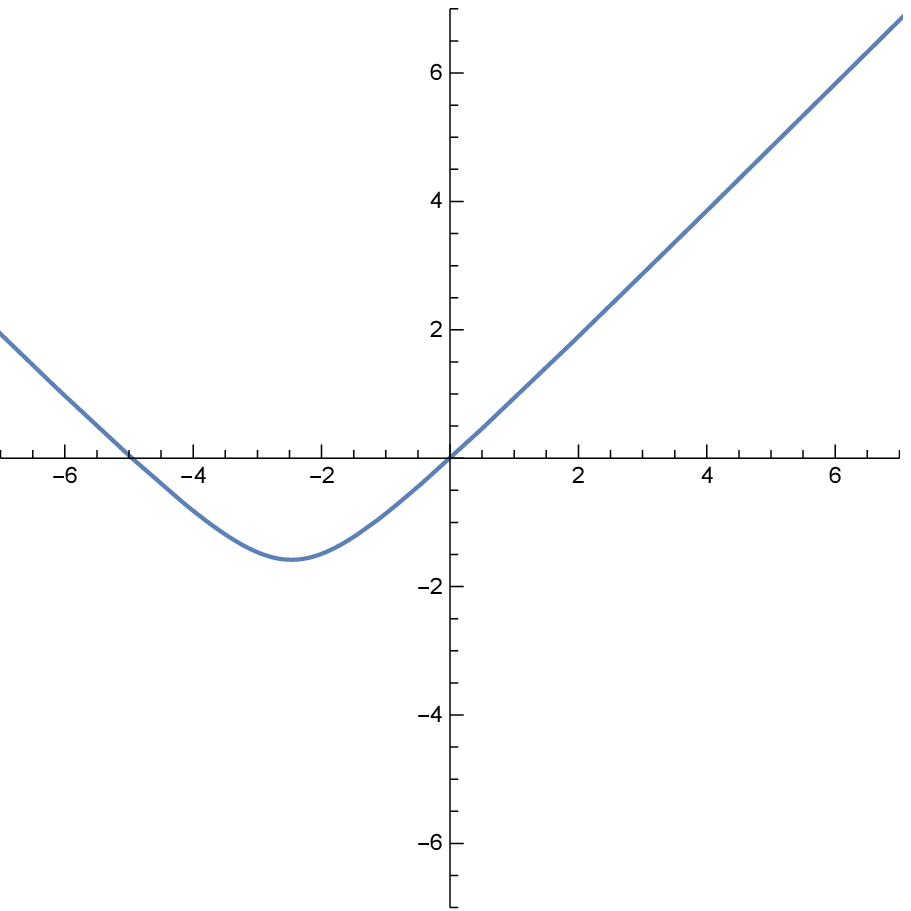}
 }
 \centering
 \subfigure[]{
 \includegraphics[width=0.3\linewidth]{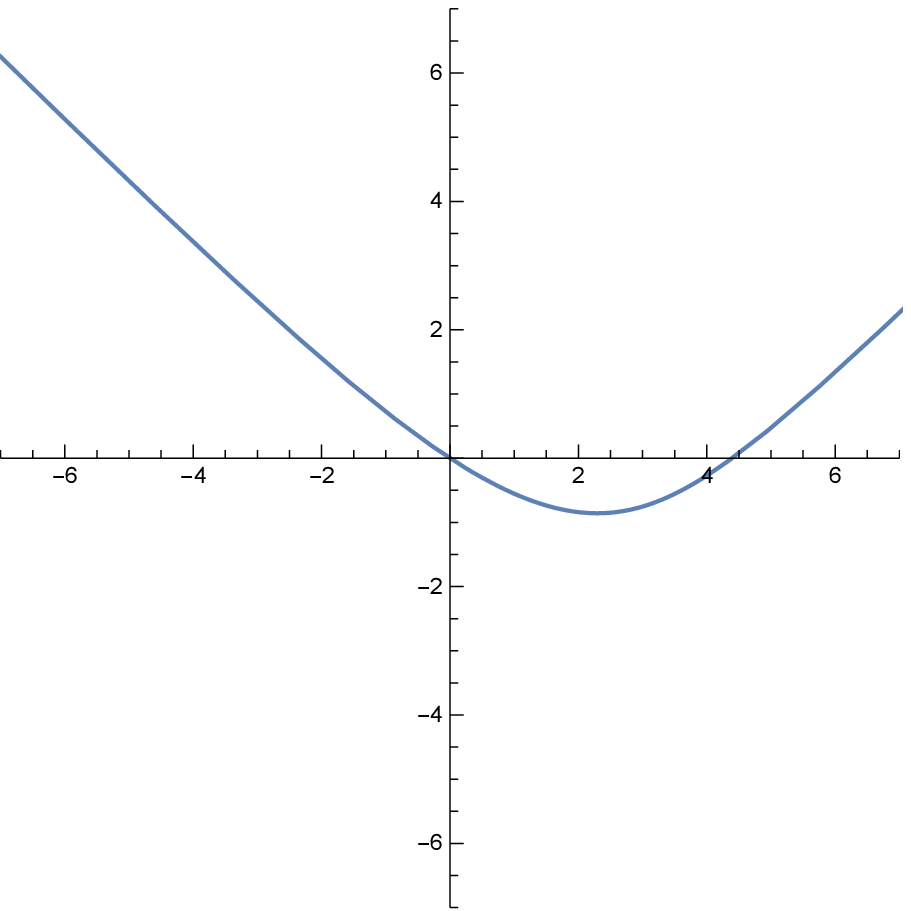}
 }
 \centering
 \subfigure[]{
 \includegraphics[width=0.3\linewidth]{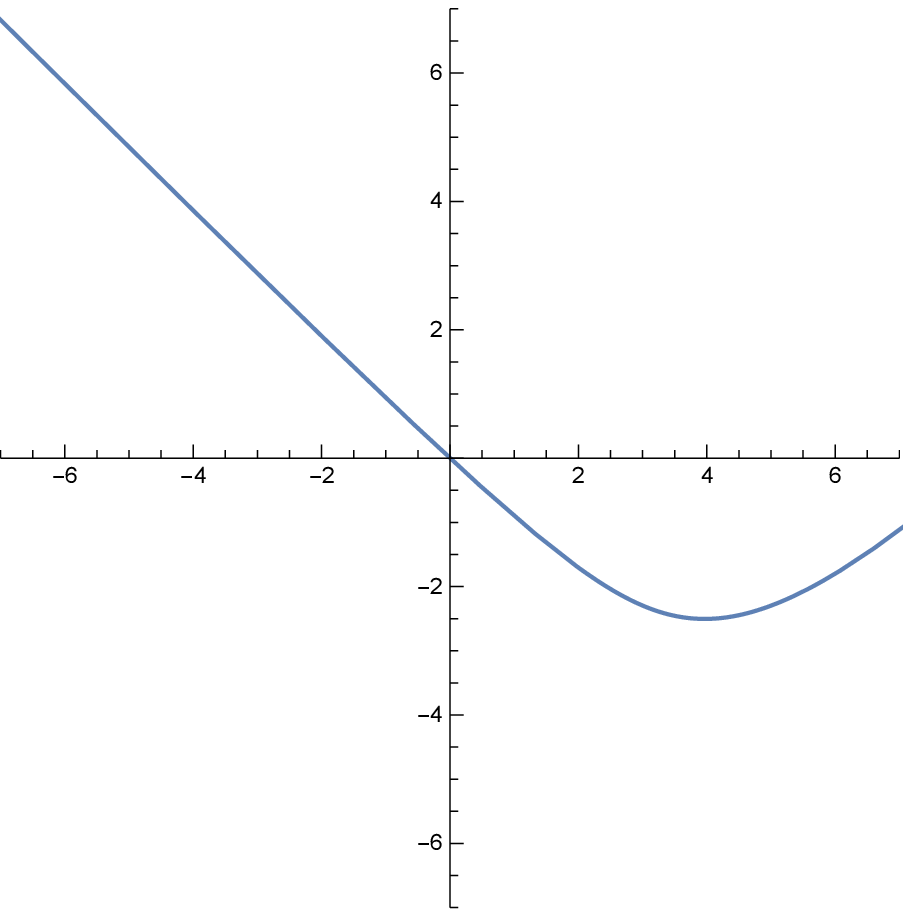}
 }
 \caption{Profiles of the defocusing mKdV flow \eqref{mkdvflow}
 given as \eqref{representation} with $\theta$ in \eqref{theta:regular}.
 Parameters in Proposition \ref{prop:regular} are
 $N=2$, $p_1=0.5$, $p_2=0.7$, $p_3=0.9$, $\alpha_1=\alpha_2=\alpha_3=\beta_1=1$,
 $\beta_2=-1$ and $\beta_3=0$, and
 $t=-18$ (a), $t=-5$ (b), $t=8$ (c).}
 \label{figure:regular_curve}
\end{figure}
\begin{figure}[H]
 \centering
 \subfigure[]{
 \includegraphics[width=0.3\linewidth]{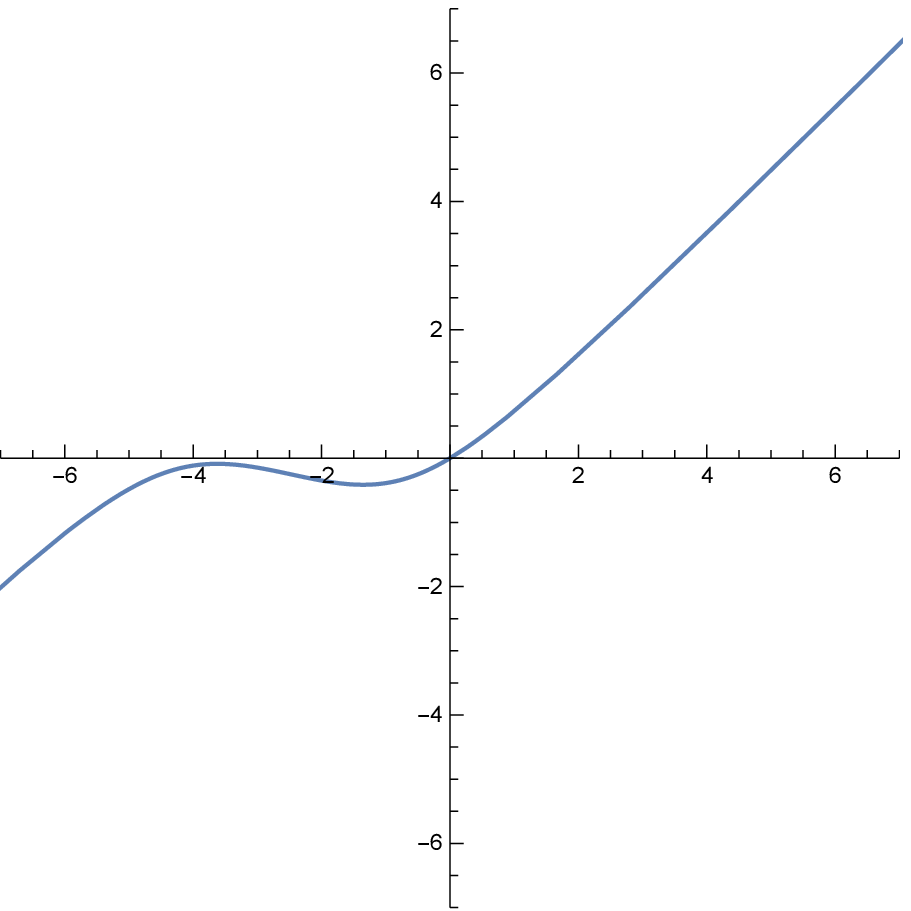}
 }
 \centering
 \subfigure[]{
 \includegraphics[width=0.3\linewidth]{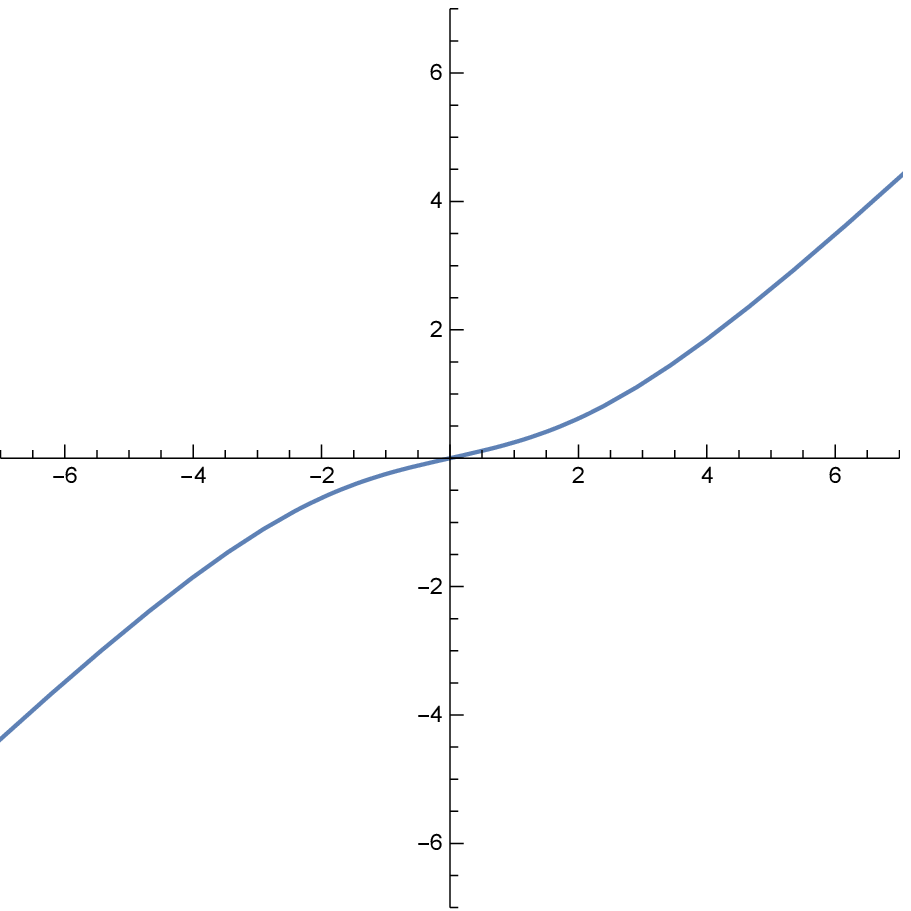}
 }
 \centering
 \subfigure[]{
 \includegraphics[width=0.3\linewidth]{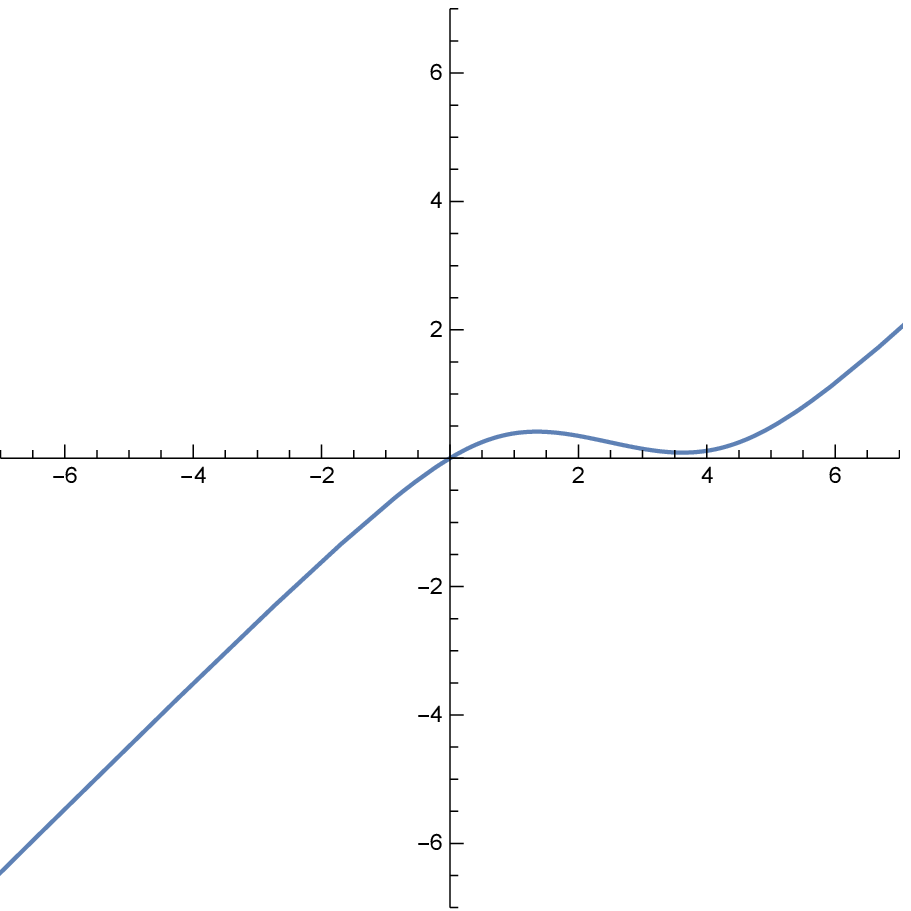}
 }
 \caption{Profiles of the defocusing mKdV flow \eqref{mkdvflow}
 given as \eqref{representation} with $\theta$ in \eqref{theta:regular}.
 Parameters in Proposition \ref{prop:regular} are
 $N=2$, $p_1=0.5$, $p_2=0.7$, $p_3=0.9$, $\alpha_1=\alpha_2=\alpha_3=\beta_1=\beta_3=1$
 and $\beta_2=-1$, and
 $t=-15$ (a), $t=0$ (b), $t=15$ (c).}
 \label{figure:regular_curve2}
\end{figure}
\section{Concluding remarks}
In this paper, we formulated the motion of spacelike curves on the Minkowki plane
preserving the arc length which is governed by the defocusing mKdV equation.
Then we constructed two classes of exact solutions to the defocusing mKdV equation
in terms of the $\tau$ functions.
Especially, one of them can be used to construct the explicit formula
for the corresponding defocusing mKdV flow in terms of the same $\tau$ functions.
However those solutions contains singular points where the $\tau$ functions have zeros.
Since we are usually interested in regular solutions,
we have also presented the regular solutions to the defocusing mKdV equation
by using different type of the $\tau$ functions and choosing suitable parameters.
These solutions describe the solitons running on a shock wave,
including the dark solitons as the special cases,
whose behavior is different from solutions of other soliton equations.
On the other hand, it seems that this class of solutions does not allow
the similar explicit formula to the solutions mentioned above.
Therefore we used a numerical integration to observe the dynamics of the
corresponding mKdV flow on the Minkowski plane curves.
It may be an interesting and important problem to extend the motion of curves
in this paper to those of discrete curves.

\section*{Acknowledgments}

This work has been partially supported by JSPS KAKENHI Grant Numbers JP16H03941, JP18H01130, JP17H06127,
JP15K04834, JP18K03435, JP15K04909, JP16K13763 and JP15K04862, and
by JST CREST Grant Number JPMJCR14D4.
One of the authors (H.P) acknowledges the support from the
``Leading Program in Mathematics for Key Technologies''
of Kyushu University.

%
%


\begin{thebibliography}{99}
  \bibitem{hasimoto:1972}
  H. Hasimoto,
  A soliton on a vortex filament,
  {\it J. Fluid Mech.} {\bf 51} (1972), 477--485.

  \bibitem{goldstein:1991}
  R. E. Goldstein and D. M. Petrich,
  The Korteweg-de Vries hierarchy as dynamics of closed curves in the plane,
  {\it Phys. Rev. Lett.} {\bf 67} (1991), 3203--3206.

  \bibitem{lamb:1976}
  G. L. Lamb, Jr.,
  Solitons and the motion of helical curves,
  {\it Phys. Rev. Lett.} {\bf 37} (1976), 235--237.

  \bibitem{chou:2002}
  K. -S. Chou and C. -Z. Qu,
  Integrable equations arising from motions of plane curves,
  {\it Phys. D} {\bf 162} (2002), 9--33.

  \bibitem{chou:2003}
  K. -S. Chou and C. -Z. Qu,
  Integrable equations arising from motions of plane curves. \rom{2},
  {\it J. Nonlinear Sci.} {\bf 13} (2003), 487--517.

  \bibitem{chou:2004}
  K. -S. Chou and C. -Z. Qu,
  Motions of curves in similarity geometries and Burgers-mKdV hierarchies,
  {\it Chaos Solitons Fractals} {\bf 19} (2004), 47--53.

  \bibitem{fujioka:2010}
  A. Fujioka and T. Kurose,
  Hamiltonian formalism for the higher KdV flows on the space of closed complex equicentroaffine curves,
  {\it Int. J. Geom. Methods Mod. Phys.} {\bf 7} (2010), 165--175.

  \bibitem{kajiwara:2016}
  K. Kajiwara, T. Kuroda and N. Matsuura,
  Isogonal deformation of discrete plane curves and discrete Burgers hierarchy,
  {\it Pac. J. Math. Ind.} {\bf 8} (2016), 14.

  \bibitem{pinkall:1995}
  U. Pinkall,
  Hamiltonian flows on the space of star-shaped curves,
  {\it Results Math.} {\bf 27} (1995), 328--332.

  \bibitem{hirose:}
  S. Hirose, J. Inoguchi, K. Kajiwara, N. Matsuura and Y. Ohta,
  Discrete local induction equation,
  arXiv:1708.01704 (2017).

  \bibitem{inoguchi:2012}
  J. Inoguchi, K. Kajiwara, N. Matsuura and Y. Ohta,
  Explicit solutions to the semi-discrete modified KdV equation and motion of discrete plane curves,
  {\it J. phys. A} {\bf 45} (2012), 045206.

  \bibitem{inoguchi:2012:2}
  J. Inoguchi, K. Kajiwara, N. Matsuura and Y. Ohta,
  Motion and B\"acklund transformations of discrete plane curves,
  {\it Kyushu J. Math.} {\bf 66} (2012), 303--324.

  \bibitem{park:2018}
  H. Park, K. Kajiwara, T. Kurose and N. Matsuura,
  Defocusing mKdV flow on centroaffine plane curves,
  submitted to JSIAM Letters.

  \bibitem{perelman:1974}
  T. L. Perelman, A. K. Fridman and M. M. El'yashevich,
  Modified korteweg-de vries equation in electrohydrodynamics,
  {\it Sov. Phys. JETP} {\bf 39} (1974), 643--646.

  \bibitem{perelman:1974:2}
  T. L. Perelman, A. K. Fridman and M. M. El'yashevich,
  On the relationship between the n-soliton solution of the modified korteweg-de vries equation and the kdv equation solution,
  {\it Phys. Lett.} {\bf 47A} (1974), 321--323.

  \bibitem{ohta:1989}
  Y. Ohta,
  Wronskian solutions to soliton equations,
  RIMS Kokyuroku 684 (1989) 1--17.

  \bibitem{hirota:2004}
  R. Hirota,
  {\it The direct method in soliton theory}
  (Cambridge University Press, Cambridge, 2004).

  \bibitem{jimbo:1983}
  M. Jimbo and T. Miwa,
  Solitons and infinite-dimensional Lie algebras,
  {\it Publ. Res. Inst. Math. Sci.} {\bf 19} (1983), 943--1001.

  \bibitem{maruno:1998}
  K. Maruno, K. Kajiwara and M. Oikawa,
  Casorati determinant solution for the discrete-time relativistic Toda lattice equation,
  {\it Phys. Lett. A} {\bf 241} (1998), 335--343.

  \bibitem{maruno:2006}
  K. Maruno and Y. Ohta,
  Casorati determinant form of dark soliton solutions of the discrete nonlinear Schr\"odinger equation,
  {\it J. Phys. Soc. Japan} {\bf 75} (2006), 054002.

  \bibitem{ohta:1993}
  Y. Ohta, R. Hirota, S. Tsujimoto and T. Imai,
  Casorati and discrete Gram type determinant representations of solutions to the discrete KP hierarchy,
  {\it J. Phys. Soc. Japan} {\bf 62} (1993), 1872--1886.

  \bibitem{ohta:1993:2}
  Y. Ohta, K. Kajiwara, J. Matsukidaira and J. Satsuma,
  Casorati determinant solution for the relativistic Toda lattice equation,
  {\it J. Math. Phys.} {\bf 34} (1993), 5190--5204.

  \bibitem{tsujimoto:2002}
  S. Tsujimoto,
  On a discrete analogue of the two-dimensional Toda lattice hierarchy,
  {\it Publ. Res. Inst. Math. Sci.} {\bf 38} (2002), 113--133.

  \bibitem{ueno:1984}
  K. Ueno, and K. Takasaki,
  Toda lattice hierarchy, in: Group representations and systems of differential equations,
  (Tokyo, 1982),
  {\it Adv. Stud. Pure Math.} {\bf 4} (1984), 1--95.

  \bibitem{freeman:1983}
  N. C. Freeman and J. J. C. Nimmo,
  Soliton solutions of the Korteweg-de Vries and the Kadomtsev-Petviashvili equations:
  the Wronskian technique,
  {\it Proc. Roy. Soc. London Ser. A} {\bf 389} (1983), 319--329.

  \bibitem{gesztesy:1991}
  F. Gesztesy and W. Schweiger,
  Rational KP and mKP-solutions in Wronskian form,
  {\it Rep. Math. Phys.} {\bf 30} (1991), 205-222.

  \bibitem{hirota:1988}
  R. Hirota, Y. Ohta and J. Satsuma,
  Solutions of the Kadomtsev-Petviashvili equation and the two-dimensional Toda equations,
  {\it J. Phys. Soc. Japan} {\bf 57} (1988), 1901--1904.

  \bibitem{hirota:1988:2}
  R. Hirota, Y. Ohta and J. Satsuma,
  Wronskian structures of solutions for soliton equations,
  {\it Progr. Theoret. Phys. Suppl.} {\bf 94} (1988), 59--72.
\end{thebibliography}
\end{document}